# THE $U(1)$-TOPOLOGICAL ELLIPTIC GENUS IS SURJECTIVE

TILMAN BAUER AND MAYUKO YAMASHITA

Abstract. We show that the topological elliptic genus from the cobordism ring of SU-manifolds to topological Jacobi forms lifts to connective topological Jacobi forms, and that this lift is surjective in homotopy.

Contents

*Acknowledgments*: The authors thank David Gepner and Lennart Meier for sharing their ongoing work. They also thank Ying-Hsuan Lin, Kazushi Ueda and Akihito Nakamura for helpful discussions. The work of MY at Perimeter Institute is supported in part by the Government of Canada through the Department of Innovation, Science and Economic Development and by the Province of Ontario through the Ministry of Colleges and Universities. She is also supported by Grant-in-Aid for JSPS KAKENHI Grant Number 20K14307, JST CREST program JPMJCR18T6, and the Simons Collaborations on Global Categorical Symmetries. Part of the collaboration was also supported by Fundamental Quantum Science Program in RIKEN Wako, Saitama.

## 1. Introduction

In [BM25], the first author and Meier introduced the graded ring spectrum of *topological Jacobi forms*, $(\mathrm{TJF}_m)_{m\geq 0}$, along with a stable version

$$\mathrm{TJF} = \operatorname*{colim}_{m\to\infty} \mathrm{TJF}_m,$$

where the latter is an $E_\infty$-algebra over the periodic spectrum TMF of topological modular forms.

Much as TMF is a spectral version of the ring MF of integral weakly holomorphic modular forms, TJF is a spectral version of the ring JF of integral weakly holomorphic Jacobi forms (see Definition 2.11 for our convention) as defined in [EZ85]. Also, in analogy with the connective version tmf, which is the spectral version of the ring of holomorphic modular forms, there is a tmf-module spectrum tjF, which is the spectral version of the ring jF of integral weak Jacobi forms.

In [LY24], the second author and Lin constructed orientations

$$\mathrm{Jac}_m\colon \mathrm{MTSU}(m) \to \mathrm{TJF}_m,$$

called *topological elliptic genera,*[1] refining the classical two-variable elliptic genus

$$\mathrm{jac}_{\mathrm{clas}}\colon \mathrm{MSU}_* \to \mathrm{jF}.$$

Here $\mathrm{MTSU}(m)$ denotes the bordism spectrum of smooth manifolds with stable tangential $\mathrm{SU}(m)$-structures. Our first result is a lift of $\mathrm{Jac}_m$ to the connective spectrum $\mathrm{tjF}_m$:

**Theorem 1.1.** *The elliptic genera* $\mathrm{Jac}_m$ *factor through orientations*

$$\mathrm{jac}_m\colon \mathrm{MTSU}(m) \to \mathrm{tjF}_m.$$

*that are compatible with the canonical inclusions* $\mathrm{MTSU}(m) \to \mathrm{MTSU}(m+1)$ *and* $\mathrm{tjF}_m \to \mathrm{tjF}_{m+1}$.

We call the stabilization

$$\mathrm{jac}\colon \mathrm{MTSU} \simeq \mathrm{MSU} \to \mathrm{tjF},$$

the *connective topological elliptic genus.*

It is known that not every weak Jacobi form occurs as an elliptic genus. For instance, the elliptic genus of a K3 surface, a Jacobi form of index 1 and weight 4, is divisible by 2, but its half is not the genus of any SU-manifold.

Our main result is:

**Theorem 1.2.** *The connective topological elliptic genus induces a surjective map in homotopy:*

$$\mathrm{jac}\colon \mathrm{MSU}_* \twoheadrightarrow \pi_*\mathrm{tjF}$$

In classical terms, it is known and quite elementary [Tot00] that

$$\mathrm{jac}_{\mathrm{clas}}\colon \mathrm{MSU}_*\left[\tfrac{1}{2}\right] \to \mathrm{jF}\left[\tfrac{1}{2}\right]$$

is surjective. However, this surjectivity fails to be true without inverting 2. Theorem 1.2 gives an exact description of the divisibility of elliptic genera as the index of the image of $\pi_*\mathrm{tjF}$ in $\mathrm{jF}_*$. Concretely, combining our main result of this paper and the computation of $\pi_*\mathrm{tjF}$ in [BM25] (Corollary 5.8), we immediately obtain an identification of the image of the classical elliptic genera for SU-manifolds.

[1] In [LY24], it is called $U(1)_m$-topological elliptic genus and denoted by $\mathrm{Jac}_{U(1)_m}$.

**Corollary 1.3.** *Let* $e_{\mathrm{jF}}\colon \pi_*\mathrm{tjF} \to \mathrm{jF}$ *denote the edge map of the descent spectral sequence. Then*

$$\mathrm{Im}\,(\mathrm{jac}_{\mathrm{clas}}\colon \mathrm{MSU}_* \to \mathrm{jF}) = \mathrm{Im}\,(e_{\mathrm{jF}}\colon \pi_*\mathrm{tjF} \to \mathrm{jF})\,.$$

*The right hand side has a concrete description: it coincides with the subring of*

$$\mathrm{jF} \cong \mathbb{Z}[b_2, b_3, b_4, b_8]/(4b_8 + b_4^2 - b_2b_3^2)$$

*generated by* $2b_2$*,* $b_3$*,* $b_4$*,* $b_2^2$*,* $b_2b_3$*,* $b_2b_4$*, and* $b_8$*.*

*Remark* 1.4. The cokernel of the map $e_{\mathrm{jF}}$ is given by

$$\mathrm{Coker}\,(e_{\mathrm{jF}}\colon \pi_*\mathrm{tjF} \to \mathrm{jF}) = \mathbb{Z}/2 \cdot \left\{b_2^{2m+1}b_8^n\right\}_{m,n\in\mathbb{Z}_{\geq 0}}.$$

⌟

Our result is in the spirit, but by no means a direct consequence, of the surjectivity in homotopy of the string orientation

$$\sigma\colon \mathrm{MString} \to \mathrm{tmf}$$

which was claimed by Hopkins and Mahowald and proven by Devalapurkar [Dev19].

We also study the kernel of jac. After inverting 2, the kernel was completely identified by Totaro [Tot00]: he defines an ideal $I \subset \mathrm{MSU}_*$ generated by "SU-flops", and shows that the classical elliptic genus gives an isomorphism

$$\mathrm{jac}_{\mathrm{clas}} \otimes \mathbb{Z}[\tfrac{1}{2}]\colon (\mathrm{MSU}_*/I)[\tfrac{1}{2}] \simeq \mathrm{jF}[\tfrac{1}{2}].$$

We make steps toward an integral refinement of Totaro's result, with the classical elliptic genus replaced by the topological elliptic genus. We show that jac is invariant under homotopy theoretically defined SU-flops so that $I \subset \ker(\mathrm{jac})$ (Theorem 8.6). This, combined with Theorem 1.2, implies that jac induces a well-defined surjection

$$\mathrm{jac}\colon \mathrm{MSU}_*/I \twoheadrightarrow \pi_*\mathrm{tjF}.$$

We do not know if this map is injective or if a modification is needed to produce an isomorphism $\mathrm{MSU}_*/I \cong \pi_*\mathrm{tjF}$.

We note that the essence of all of these results is 2-local: all of the torsion in MSU and tjF is 2-torsion, and away from 2, all these homology theories are complex orientable and evenly graded. This is not true unstably (i.e. for $\mathrm{MTSU}(m)$ and $\mathrm{tjF}_m$), but the unstable elliptic genera are not surjective in this cases even at $p > 2$.

The paper is organized as follows. Section 2 surveys the construction of TJF and Section 3 the construction of the (nonconnective) topological elliptic genera and their algebraic counterparts. In Section 4, we construct the connective topological elliptic genera $\mathrm{jac}_m$. In Section 5, we study the descent (or Adams–Novikov) spectral sequence and the resulting homotopy groups of tjF. After the preparatory Section 6 on MSU, we prove the main result in Section 7. Section 8 is on the SU-flop invariance of $\mathrm{jac}_m$.

## 2. Topological Jacobi forms

This chapter contains a brief review of the construction of topological Jacobi forms. More details can be found in [GM23, BM25].

Topological Jacobi forms are a graded TMF-algebra spectrum $(\mathrm{TJF}_m)_{m\geq 0}$ in the sense that there is a TMF-bilinear multiplication

$$\mathrm{TJF}_{m_1} \otimes_{\mathrm{TMF}} \mathrm{TJF}_{m_2} \to \mathrm{TJF}_{m_1+m_2}$$

which is $E_2$, in particular homotopy commutative. These spectra are constructed using equivariant TMF. In [GM23, GLP24], the authors construct a global equivariant spectrum $\{\mathrm{TMF}_G\}$, where $G$ ranges through all compact Lie groups. We will be interested in the case $G = U(1)$. For any $G$ and (real) $G$-representation $V$, we denote by

$$\mathrm{TMF}[V]^G = (\mathrm{TMF}_G \otimes S^V)^G$$

the (genuine) $G$-fixed points of the smash product of $\mathrm{TMF}_G$ and the one-point compactification of $V$. This definition extends to virtual $G$-representation by defining $S^{W-V} = S^W \otimes DS^V$, where $D$ denotes the Spanier-Whitehead dual. For $G = U(1)$, let $L$ denote the standard 1-dimensional complex representation.

**Definition 2.1.** Let $m$ be any integer. Then $\mathrm{TJF}_m$ is the TMF-module spectrum

$$\mathrm{TJF}_m := \mathrm{TMF}[mL]^{U(1)}.$$

There is a more (spectral) algebro-geometric description of these cohomology theories, which will be important in constructing descent spectral sequences. This is using the language of spectral algebraic geometry [Lur18c, Lur18a, Lur18b]. The reader unfamiliar with the details of this may wish to take the descent spectral sequence for granted; the rest of this paper is formulated in the language of classical (equivariant) homotopy theory.

In [GM], Gepner and Meier construct for a given oriented spectral elliptic curve $E$ over a spectral base scheme $S$ a functor

$$\mathrm{Ell}\colon \mathrm{Spc}_{\mathrm{Orb}} \to \mathrm{Top}^{\mathrm{loc}}_{\mathrm{CAlg}}$$

from the $\infty$-category of orbispaces to the $\infty$-category of locally spectrally ringed topological spaces with the following properties:

- $\mathrm{Ell}(\mathbb{B}U(1)^n) \cong (E, \mathcal{O}_E)^{\times_S n}$, where $\mathbb{B}G = [*/G] \in \mathrm{Spc}_{\mathrm{Orb}}$ is the global stack quotient of $G$ acting trivially on a point, for a compact Lie group $G$.
- Ell commutes with colimits;
- $\mathrm{Ell} = \mathrm{Ell}_S$ satisfies étale descent for finite global quotients, i.e. given an étale map $S' \to S$ and a finite $G$-complex $X$, the map

$$\mathrm{Ell}_{S'}([X/G]) \cong \mathrm{Ell}_S([X/G]) \times_S S' \to \mathrm{Ell}_S([X/G])$$

  is étale as well. Here again, $[X/G] \in \mathrm{Spc}_{\mathrm{Orb}}$ denotes the stack quotient of $X$ by $G$.

The last property allows us to extend $\mathrm{Ell}_S$ to bases that are (nonconnective) spectral Deligne–Mumford stacks such as the stack $S = \mathcal{M}^{\mathrm{or}}_{\mathrm{ell}}$ of all oriented spectral elliptic curves. For the remainder of this paper, we will exclusively consider this base. The resulting functor will take values in spectral DM–stacks over $\mathcal{M}^{\mathrm{or}}_{\mathrm{ell}}$ that are relative schemes, at least when applied to finite global quotients.

*Remark* 2.2. In [GM23], the authors use a different target category for the functor Ell, resulting in a different version of Ell for nonabelian groups of invariance. There is, however, no difference after applying the functor of global sections, $(X, \mathcal{O}_X) \mapsto \Gamma_X(\mathcal{O}_X)$. ⌟

The stack $\mathrm{Ell}(\mathbb{B}U(1))$ is the universal oriented spectral elliptic curve $\mathcal{E}^{\mathrm{or}}$. The canonical map $p\colon [X/G] \to \mathbb{B}G$ gives rise to a functor

$$\mathcal{E}\mathrm{ll}_G\colon \mathrm{Spc}_G \xrightarrow{X\mapsto[X/G]} \mathrm{Spc}_{\mathrm{Orb}\,/\mathbb{B}G} \xrightarrow{\mathrm{Ell}} \mathrm{Shv}_{\mathcal{M}^{\mathrm{or}}_{\mathrm{ell}}/\mathrm{Ell}(\mathbb{B}G)} \xrightarrow{(M,\mathcal{O}_M)\mapsto p_*\mathcal{O}_M} \mathrm{QCoh}(\mathrm{Ell}(\mathbb{B}G))^{\mathrm{op}}.$$

The equivariant cohomology theory $\mathrm{TMF}_G$ is obtained by composing with the pushforward along $q\colon \mathbb{B}G \to *$,

$$\mathrm{Spc}_G \xrightarrow{\mathcal{E}\mathrm{ll}_G} \mathrm{QCoh}(\mathrm{Ell}(\mathbb{B}G))^{\mathrm{op}} \xrightarrow{q_*} \mathrm{QCoh}(\mathcal{M}^{\mathrm{or}}_{\mathrm{ell}})^{\mathrm{op}} \xrightarrow{\Gamma} \mathrm{Mod}^{\mathrm{op}}_{\mathrm{TMF}}.$$

For $G = U(1)$, we define the quasicoherent sheaves

$$\mathcal{O}_{\mathcal{E}^{\mathrm{or}}}(-ke) = \mathcal{E}\mathrm{ll}(S^{kL}) \in \mathrm{QCoh}(\mathcal{E}^{\mathrm{or}})$$

By [GM23, Lemma 8.1], these sheaves are invertible line bundles, allowing us to extend this to negative $k$. Homotopically, this means that $\mathcal{E}\mathrm{ll}$ extends to $G$-spectra and is symmetric monoidal, and the inverse line bundles are given by applying Spanier-Whitehead duality to the representation spheres.

**Proposition 2.3.**

$$\mathrm{TJF}_m \simeq \Gamma(\mathcal{E}^{\mathrm{or}}, \mathcal{O}_{\mathcal{E}^{\mathrm{or}}}(me)) \simeq \Gamma(\mathcal{M}_{\mathrm{ell}}, q_*\mathcal{O}_{\mathcal{E}^{\mathrm{or}}}(me))$$

We have the following additive description of $\mathrm{TJF}_m$ [BM25]:

**Theorem 2.4.** *There is an equivalence of* TMF*-module spectra*

$$\mathrm{TJF}_m \simeq \mathrm{TMF} \otimes P_m,$$

*where* $P_m = \mathrm{cofib}(\widetilde{\mathrm{tr}}\colon \Sigma\mathbb{CP}^{m-1} \to S^0)$ *is the cofiber of the reduced* $U(1)$*-transfer.*

Letting $L$ denote the tautological complex line bundle on $\mathbb{CP}^m$, the cofiber of the *unreduced* $U(1)$-transfer is identified as

$$\mathrm{cofib}\left(\mathrm{tr}\colon \Sigma\mathbb{CP}^{m-1}_+ \to S^0\right) \simeq \Sigma^2\mathbb{CP}^m_{-1} = (\mathbb{CP}^m)^{-L} \simeq D((\mathbb{CP}^m)^{-mL}) =_{\mathrm{def}} D(\mathbb{CP}^0_{-m}),$$

where the first isomorphism follows by [BM25, Proposition A.11], and the latter isomorphism is by Atiyah duality and the stable equivalence $T\mathbb{CP}^m \simeq (m+1)L - 1$ of complex virtual vector bundles.

Thus $P_m$ looks like the dual of $\mathbb{CP}^0_{-m}$ with the 2-dimensional cell removed. That is, $P_m$ fits into a fiber sequence

$$P_m \xrightarrow{f_m} D(\mathbb{CP}^0_{-m}) \to S^2. \tag{2.5}$$

Multiplication with the equivariant Euler class of the fundamental representation $L \in \mathrm{Rep}(U(1))$,

$$\chi(L) \in \pi_0\mathrm{TJF}_1 \tag{2.6}$$

gives the *stabilization map*,

$$a\colon \mathrm{TJF}_m \xrightarrow{x\mapsto\chi(L)x} \mathrm{TJF}_{m+1} \tag{2.7}$$

Under the isomorphism of Prop. 2.3, the stabilization map is identified with the canonical map

$$\Gamma(\mathcal{E}^{\mathrm{or}}, \mathcal{O}_{\mathcal{E}^{\mathrm{or}}}(me)) \to \Gamma(\mathcal{E}^{\mathrm{or}}, \mathcal{O}_{\mathcal{E}^{\mathrm{or}}}((m+1)e)).$$

It was also shown in [BM25] that the multiplication induced by the map $S^{mV} \otimes S^{m'V} \to S^{(m+m')V}$ equips $\bigoplus_{m\geq 0}\mathrm{TJF}_m$ with the structure of a filtered $E_2$-ring spectrum over TMF.

In this paper, we are particularly interested in the stable version of topological Jacobi forms:

**Definition 2.8** (TJF)**.** We define[2]

$$\mathrm{TJF} := \operatorname*{colim}_{m\to\infty}\mathrm{TJF}_m \simeq \operatorname*{colim}_{m\to\infty}\Gamma(\mathcal{E}^{\mathrm{or}};\mathcal{O}_{\mathcal{E}^{\mathrm{or}}}(me)) \simeq \Gamma(\mathcal{E}^{\mathrm{or}}-\{e\};\mathcal{O}_{\mathcal{E}^{\mathrm{or}}-\{e\}}),$$

where the colimit is taken over (2.7). The equivalence with sections over the open substack $\mathcal{E}^{\mathrm{or}}-\{e\}$, the complement of the unit section, comes from the fact that a meromorphic function with a pole of unlimited order at $e$ and no other poles is the same as a holomorphic function defined away from $e$.

We also define a connective version

$$\mathrm{tjF}_m := \mathrm{tmf}\otimes P_m, \quad \mathrm{tjF} := \mathrm{tmf}\otimes P_\infty$$

Note that we are not claiming more than a tmf-module spectrum structure on tjF. They come with an equivalence of tmf-module spectra

$$\mathrm{tjF}_m[\Delta^{-24}] \simeq \mathrm{TJF}_m, \quad \mathrm{tjF}[\Delta^{-24}] \simeq \mathrm{TJF}$$

simply because $\mathrm{tmf}[\Delta^{-24}] \simeq \mathrm{TMF}$.

2.1. **Descent spectral sequences.** Given a spectral Deligne-Mumford stack $\mathcal{M}$ with a quasicoherent $\mathcal{O}_{\mathcal{M}}$-module $\mathcal{F}$, there is a conditionally convergent spectral sequence, called the *descent spectral sequence*,

$$H^p(\mathcal{M}^\heartsuit;\pi_q(\mathcal{F})) \Longrightarrow \pi_{q-p}\Gamma(\mathcal{F}), \tag{2.9}$$

where $\mathcal{M}^\heartsuit$ denotes the underlying classical DM-stack. This spectral sequence is multiplicative if $\mathcal{F}$ is a sheaf of ring spectra.

To construct it, let $U\to\mathcal{M}$ be an étale cover by a spectral scheme. By the sheaf condition,

$$\Gamma(\mathcal{M};\mathcal{F}) \simeq \mathrm{Tot}(\Gamma(U\times_{\mathcal{M}}U\times\cdots\times_{\mathcal{M}}U;p^*\mathcal{F}))$$

where $p\colon U\times_{\mathcal{M}}U\times\cdots\times_{\mathcal{M}}U\to\mathcal{M}$ is the projection map. The descent spectral sequence is the Bousfield–Kan spectral sequence associated with this cosimplicial spectrum. In particular, its $E_2$-term is the Čech cohomology

$$\check{H}^p_{U^\heartsuit}(\mathcal{M}^\heartsuit;\pi_q\mathcal{F})$$

associated to the cover $U^\heartsuit\to X^\heartsuit$.

We have two flavors of descent spectral sequences:

$$H^p(\mathcal{M}_{\mathrm{ell}};\pi_q(\mathcal{O}_{\mathcal{M}^{\mathrm{or}}_{\mathrm{ell}}})) \Longrightarrow \pi_{q-p}\mathrm{TMF} \qquad (\pi_{2q}\mathcal{O}_{\mathcal{M}^{\mathrm{or}}_{\mathrm{ell}}}\cong\omega^q)$$
$$H^p((\mathcal{E}^{\mathrm{or}})^\heartsuit-\{0\};\pi_q(\mathcal{O}_{\mathcal{E}^{\mathrm{or}}-\{0\}})) \Longrightarrow \pi_{q-p}\mathrm{TJF}$$

Each of these is multiplicative and the descent spectral sequence for TJF is a module spectral sequence over the one for TMF. The stack $\mathcal{E}=(\mathcal{E}^{\mathrm{or}})^\heartsuit$ is the universal elliptic curve, and $\pi_{2q}(\mathcal{O}_{\mathcal{E}^{\mathrm{or}}})$ is the pullback of $\omega^q$ to $\mathcal{E}$.

---

[2] In [BM25] and [LY24], the notation $\mathrm{TJF}_\infty$ is used for this object. We omit the $\infty$ in the subscript.

If $\mathrm{pr}\colon \mathcal{E}^{\mathrm{or}} \to \mathcal{M}^{\mathrm{or}}_{\mathrm{ell}}$ denotes the projection, we trivially have that $\Gamma_{\mathcal{E}^{\mathrm{or}}}(\mathcal{O}_{\mathcal{E}^{\mathrm{or}}}) \cong \Gamma_{\mathcal{M}_{\mathrm{ell}}}(\mathrm{pr}_*\mathcal{O}_{\mathcal{E}^{\mathrm{or}}})$. Since $\mathcal{E}^{\mathrm{or}} - \{0\}$ is fiberwise affine, $\mathrm{pr}_*$ is an exact functor and we obtain an identification of $E_2$-terms

$$H^p((\mathcal{E}^{\mathrm{or}})^{\heartsuit} - \{0\}; \pi_q(\mathcal{O}_{\mathcal{E}^{\mathrm{or}}-\{0\}})) \cong H^p(\mathcal{M}_{\mathrm{ell}}; \pi_q(\mathrm{pr}_*\mathcal{O}_{\mathcal{E}^{\mathrm{or}}-\{0\}})).$$

Locally, over a Weierstrass curve

$$E(a_1, \cdots, a_6, x, y)\colon y^2 + a_1xy + a_3y = x^3 + a_2x^2 + a_4x + a_6 \tag{2.10}$$

over the ring $A = \mathbb{Z}[a_1, a_2, a_3, a_4, a_6, \Delta^{-1}]$, the module sheaf $\pi_*(\mathrm{pr}_*\mathcal{O}_{\mathcal{E}^{\mathrm{or}}-\{0\}})$ corresponds to the module $A[x, y]/E(a_1, \dots, a_6, x, y)$.

2.2. **The relation with classical Jacobi forms.** Recall the definitions of Jacobi forms from [EZ85]:

**Definition 2.11** (Jacobi forms)**.** We denote by $\mathfrak{H} := \{\tau \in \mathbb{C} \mid \mathrm{Im}(\tau) > 0\}$ the upper half space of the complex plane. For each $m \in \mathbb{Z}_{\geq 0}$ and $k \in \mathbb{Z}$, consider holomorphic functions of $(z, \tau) \in \mathbb{C} \times \mathfrak{H}$ satisfying the transformation properties

$$\phi\left(\frac{a\tau + b}{c\tau + d}, \frac{z}{c\tau + d}\right) = (c\tau + d)^k e^{\frac{\pi i m c z^2}{c\tau + d}} \phi(\tau, z),$$
$$\phi(\tau, z + \lambda\tau + \mu) = e^{-\pi i m(\lambda^2\tau + 2\lambda z + \lambda + \mu)}\phi(\tau, z)$$

for all $\left(\begin{smallmatrix} a & b \\ c & d \end{smallmatrix}\right) \in SL(2, \mathbb{Z})$ and $(\lambda, \mu) \in \mathbb{Z}^2$, and having Fourier expansions

$$\phi(q, y) = \sum_{r \in \mathbb{Z} + \frac{m}{2}} \sum_{n \geq N} c(n, r) q^n y^r$$

where $(q, y) = (\exp(2\pi i\tau), \exp(2\pi i z))$ for some integer $N$. Such functions are called *weakly holomorphic* Jacobi forms of index $\frac{m}{2}$ and weight $k$. If $c(n, r) \neq 0$ only when $n \geq 0$, then such functions are called *weak* Jacobi forms.[3] If all $c(n, r) \in \mathbb{Z}$, we add the adjective *integral* in all the above cases.

We denote by $\mathrm{JF}_m$ (and $\mathrm{jF}_m$) the set of all weakly holomorphic (resp. weak) integral Jacobi forms of index $\frac{m}{2}$. We equip $\mathrm{JF}_m$ and $\mathrm{jF}_m$ with a $\mathbb{Z}$-grading ("degree") whose degree-$n$ part consists of Jacobi forms with weight $k = \frac{n}{2} - m$. This convention makes $\mathrm{jF}_m$ a $\mathbb{Z}$-graded module over the $\mathbb{Z}$-graded ring mf of integral holomorphic modular forms (with grading given by weight/2), and likewise $\mathrm{JF}_m$ over the ring MF of integral weakly-holomorphic modular forms.

The following was proved in [BM25, Thm 2.5]:

**Lemma 2.12.** *We have*

$$\mathrm{JF}_m \cong H^0((\mathcal{E}^{\mathrm{or}})^{\heartsuit}; \pi_*\mathcal{O}_{\mathcal{E}^{\mathrm{or}}}(me)).$$

□

Topological Jacobi Forms are a homotopical refinement of the classical ring of weakly holomorphic integral Jacobi Forms, in the sense that they are equipped with an edge map coming from the descent spectral sequence (2.9),

$$e_{\mathrm{JF}}\colon \pi_*\mathrm{TJF}_m \to \mathrm{JF}_m$$

This map is a special case of the map from $G$-equivariant TMF to $G$-equivariant weakly holomorphic modular forms, which are referred to as *character maps* in [LY24].

[3] There is also a notion of *holomorphic Jacobi forms*, but we do not use it in this paper.

We have multiplication maps $\mathrm{jF}_k \otimes_{\mathrm{mf}} \mathrm{jF}_m \to \mathrm{jF}_{k+m}$, which makes $\mathrm{jF}_\bullet = \bigoplus_m \mathrm{jF}_m$ into a $\mathbb{Z}_{\geq 0}$-graded mf-algebra. In [GW20, Theorem 2.7], Gritsenko gave a complete computation of $\mathrm{jF}_\bullet$ as an mf-algebra. In particular, the Euler class $\chi(L)$ (2.6) maps via the character map $e_{\mathrm{JF}}$ to a class $a \in (\mathrm{jF}_1)$ of degree 0. This class is called $\phi_{-1,\frac{1}{2}}$ in [GW20] and can be expressed in terms of classical Jacobi theta functios and the Dedekind eta function as

$$a = \phi_{-1,\frac{1}{2}} = \frac{\theta_{11}(z,\tau)}{\eta^3(\tau)}.$$

In analogy with the stabilization maps (2.7), we thus get the stabilization maps of Jacobi forms,

$$a\cdot \colon \mathrm{jF}_m \to \mathrm{jF}_{m+1}$$

**Definition 2.13.** Define the stabilized ring of weak Jacobi forms to be

$$\mathrm{jF} = \mathrm{jF}_\infty := \operatorname*{colim}_{m\to\infty} \left( \mathrm{jF}_m \xrightarrow{a\cdot} \mathrm{jF}_{m+1} \xrightarrow{a\cdot} \cdots \right).$$

It follows from [GW20] that

$$\mathrm{jF} = \mathbb{Z}[b_2, b_3, b_4, b_8]/(4b_8 + b_4^2 - b_2 b_3^2).$$

An independent computation of this ring stems from [BM25], where also the higher cohomology of the universal elliptic curve is computed. For details, see Section 5.1.

## 3. Topological Elliptic Genera

In this section, we review some facts about the topological elliptic genera introduced in [LY24] and a different version introduced in a somewhat different setting by Ando–French–Ganter [AFG08]. We refer the reader to those papers for details. We also study multiplicativity properties not addressed in the cited papers.

### 3.1. **Tangential SU-bordism.**

For any representation $V$ of a Lie group $G$ over $k$, let $\overline{V}$ be the reduced virtual 0-dimensional representation $V - \dim_k V$ in the representation ring.

Denote by $\mathrm{MTSU}(m)$ the *tangential* $\mathrm{SU}(m)$-bordism spectrum, defined to be the Thom spectrum

$$\mathrm{MTSU}(m) := B\mathrm{SU}(m)^{-\overline{V}}.$$

of the negative of the reduced universal bundle $\overline{V}$ over $B\mathrm{SU}(m)$. Its homotopy groups are the bordism groups of closed manifolds equipped with stable $\mathrm{SU}(m)$-structures on their *tangent* bundles.

The collection of $\mathrm{MTSU}(m)$ becomes a graded ring spectrum by virtue of the multiplication maps

$$\mathrm{MTSU}(m_1) \otimes \mathrm{MTSU}(m_2) \to \mathrm{MTSU}(m_1 + m_2),$$

representing the cartesian product of manifolds, as well as the stabilization maps

$$a\colon \mathrm{MTSU}(m) \to \mathrm{MTSU}(m+1)$$

induced by the inclusion $\mathrm{SU}(m) \hookrightarrow \mathrm{SU}(m+1)$. Inverting $a$, we recover the more familiar bordism spectrum of *normal* SU-manifolds,

$$\mathrm{MSU} \simeq \mathrm{MTSU} := \operatorname*{colim}_{m\to\infty} \mathrm{MTSU}(m).$$

In fact, $\bigoplus_{m\geq 0} \mathrm{MTSU}(m)$ has the structure of a filtered $E_1$-ring spectrum. To see this, it is convenient to consider the $\infty$-categorical point of view of [ABG10] for constructing Thom spectra. Let $X$ be an $\infty$-groupoid (i.e. a space) with a functor $F$ to the category $\mathrm{Pic}_S$ of invertible spectra (i.e. spheres). Its colimit is the (generalized) Thom spectrum $X^F$. If $X$ is an $E_n$-monoidal $\infty$-groupoid (i.e. an $E_n$-space) and $F$ is $E_n$-monoidal, then $X^F$ inherits an $E_n$-structure. Moreover, an $E_n$-monoidal natural transformation $F \to G$ of functors from $X$ to $\mathrm{Pic}_S$ gives rise to an $E_n$-map $X^F \to X^G$ in a natural way.

In our application, $X$ is the $\infty$-groupoid associated to the monoidal topological groupoid of finite-dimensional hermitian vector spaces with a trivialization of their determinant line, so that $|X| \simeq \coprod_{m\geq 0} B\mathrm{SU}(m)$, and the monoidal functor $F$ is given by sending such a vector space $V$ to $\Sigma^{\dim_{\mathbb{R}} V} \otimes DS^V$, with $D$ being the (monoidal) Spanier-Whitehead dual. Since $F(V) = J(\dim V - V)$, the resulting $E_1$-Thom spectrum is

$$\operatorname{colim} F \simeq \bigoplus_{m\geq 0} \mathrm{MTSU}(m).$$

In the same way, $\bigoplus_m M\mathrm{SU}(m)$ arises by choosing the Spanier–Whitehead dual of the above, $F(V) = \Sigma^{-\dim_{\mathbb{R}} V} S^V$. Note that in the classical literature, what we call $M\mathrm{SU}(m)$ here would be called $\Sigma^{-2m} M\mathrm{SU}(m)$.

*Remark* 3.1. The functor $F$ above is indeed only $E_1$-monoidal. However, the functor $F(V) = DS^V$, whose colimit is the unreduced Thom spectrum, is $E_\infty$-monoidal, at least when restricting to even-dimensional vector spaces $V$. (We thank L. Meier for the observation that the braiding does not preserve the determinant line trivialization for odd-dimensional $V$.) That is, the graded spectrum $\bigoplus_m \Sigma^{-4m} \mathrm{MTSU}(2m)$ is a graded $E_\infty$-spectrum. Similarly, the spectrum $\bigoplus_m \Sigma^{-2m} \mathrm{TJF}_m$ carries a graded $E_\infty$-structure (Meier, unpublished), so that this setup might actually be more natural to consider. ⌟

3.2. **Complex-analytic elliptic genera.** Classically, elliptic genera (cf. [Gri99, Tot00]) are defined as graded maps

$$\mathrm{jac}_{\mathrm{clas},m} \colon \pi_* \mathrm{MTSU}(m) \to \mathrm{jF}_m \quad (m \geq 0)$$

stabilizing to a homomorphism of graded rings

$$\mathrm{jac}_{\mathrm{clas}} \colon \mathrm{MSU}_* \to \mathrm{jF}.$$

Over $\mathbb{C}$, for a closed stable SU-manifold $M$, the elliptic genus is given by the integral

$$\mathrm{jac}_{\mathrm{clas}}(M)(y,q) = \int_M \mathrm{Todd}(TM) \wedge \mathrm{Ch}\left(\overline{\mathbb{TM}}_{q,y}\right),$$

where $\overline{\mathbb{TM}}_{q,y}$ is the following element in $\mathrm{KU}^0(M)[[q, y^{-1}, y]]$:

$$\overline{\mathbb{TM}}_{q,y} := \bigotimes_{m\geq 0} \bigwedge\nolimits_{-q^m y^{-1}} \overline{T^*M} \otimes \bigotimes_{m\geq 1} \bigwedge\nolimits_{-q^m y} \overline{TM} \otimes \bigotimes_{m\geq 1} \mathrm{Sym}_{q^m} \overline{T^*M} \otimes \bigotimes_{m\geq 1} \mathrm{Sym}_{q^m} \overline{TM}.$$

Here all the tensor/exterior products are over $\mathbb{C}$, and $\overline{TM} := TM - \dim M$. We use the standard convention $\mathrm{Sym}_q V := \sum_n (\mathrm{Sym}^n V) q^n$ and $\bigwedge_q V := \sum_n (\bigwedge^n V) q^n$.

3.3. **Topological elliptic genera à la Lin–Yamashita.** The *topological elliptic genera* introduced in [LY24] are refinements of the complex-analytic genera to maps of spectra

$$\mathrm{Jac}_m \colon \mathrm{MTSU}(m) \to \mathrm{TJF}_m.$$

for each $m \in \mathbb{Z}_{\geq 0}$. The key ingredient in the construction is the canonical $BU\langle 6\rangle$-structure $\mathfrak{s}$ on the virtual complex representation

$$\overline{V} \boxtimes_{\mathbb{C}} \overline{L} = \operatorname*{colim}_m (V_m - m) \boxtimes_{\mathbb{C}} (L-1)$$

of $SU(\infty) \times U(1)$. Let us explain it in a way useful for our purposes. For each nonnegative integer $n$, let $\mathrm{ku}\langle n\rangle$ denote the $n$-connective cover of ku. The multiplication on $\mathrm{ku}\langle \bullet\rangle$ adjoins to a map of spectra

$$\mathrm{ku}\langle n\rangle \to \mathrm{map}(\mathrm{ku}\langle n'\rangle, \mathrm{ku}\langle n+n'\rangle),$$

where map is the mapping spectrum. Taking $\Omega^\infty$ gives a map

$$BU\langle n\rangle \to \mathrm{Map}(\mathrm{ku}\langle n'\rangle, \mathrm{ku}\langle n+n'\rangle) \to \mathrm{map}^{E_\infty}(BU\langle n'\rangle, BU\langle n+n'\rangle), \tag{3.2}$$

where $\mathrm{Map}(-,-) = \Omega^\infty \mathrm{map}(-,-)$ is the mapping *space* and the second map applies $\Omega^\infty$ on the inside.

**Definition 3.3.** Setting $n = 2$ and $n' = 4$ in (3.2) and noticing that $BU\langle 4\rangle = BSU$ and $BU\langle 2\rangle = BU$, let us consider the composition

$$BU(1) \xrightarrow{\iota} BU \xrightarrow{(3.2)} \mathrm{map}^{E_\infty}(BSU, BU\langle 6\rangle),$$

where $\iota\colon BU(1) \hookrightarrow BU$ classifies $\overline{L} = L - 1$. We define $\mathfrak{s}$ to be the map obtained by taking the adjoint of the above map, we get

$$\mathfrak{s}\colon BSU \times BU(1) \to BU\langle 6\rangle$$

The map $\mathfrak{s}$ preserves the $E_\infty$ structure on $BSU$ in the source and on the target $BU\langle 6\rangle$.

By construction, $\mathfrak{s}$ gives a $BU\langle 6\rangle$-structure on the virtual representation $\overline{V} \boxtimes_{\mathbb{C}} \overline{L}$ as desired. Moreover, the $E_\infty$-ness stated above implies that this $BU\langle 6\rangle$-structure behaves well under taking external direct sum of SU-representations.

Consider the following map in $\mathrm{Sp}_{U(1)}$:

$$\begin{aligned} \chi_m \colon \mathrm{MTSU}(m) &= B\mathrm{SU}(m)^{-\overline{V}_m} \\ &\xrightarrow{\chi(V_m \boxtimes_{\mathbb{C}} L)\cdot} B\mathrm{SU}(m)^{V_m \boxtimes_{\mathbb{C}} L - \overline{V}_m} \\ &= B\mathrm{SU}(m)^{\overline{V}_m \boxtimes_{\mathbb{C}} \overline{L} + mL}. \end{aligned}$$

Here, $\mathrm{MTSU}(m)$ and $B\mathrm{SU}(m)$ are trivial $U(1)$-spectra, and $V_m \boxtimes_{\mathbb{C}} L = V_{U(1)} \otimes_{\mathbb{C}} V_m$ is regarded as a $U(1)$-equivariant vector bundle over $B\mathrm{SU}(m)$. The arrow in the second row is given by the inclusion of the zero section of $V_m \boxtimes_{\mathbb{C}} L$.

**Lemma 3.4.** *The spectrum* $\bigoplus_m B\mathrm{SU}(m)^{\overline{V}_m \boxtimes_{\mathbb{C}} \overline{L} + mL}$ *is an* $U(1)$*-equivariant* $E_1$*-ring spectrum, and the maps* $\chi_m$ *assemble to an* $E_1$*-map.*

*Proof.* The spectrum $\bigoplus_m B\mathrm{SU}(m)^{\overline{V}_m \boxtimes_{\mathbb{C}} \overline{L} + mL}$ is the Thom construction of the symmetric monoidal functor $F\colon X \to \mathrm{Pic}_S$ with $X$, as in Subsection 3.1, the groupoid of finite-dimensional hermitian vector spaces $V$ with a trivialization of the determinant line, sending $V$ to the $U(1)$-spectrum $S^{V \otimes_{\mathbb{C}} L} \wedge S^{\dim_{\mathbb{C}}(V) L} \wedge DS^V$. The natural inclusion $S^0 \to S^{V\otimes_{\mathbb{C}} L}$ is $U(1)$-equivariant and monoidal, so that the induced map on Thom spectra is $E_1$. ☐

**Definition 3.5** ($\mathrm{Jac}_m$, [LY24, Proposition 4.19])**.** Using the unit map $\eta\colon S \to \mathrm{TMF} \in \mathrm{Sp}_{U(1)}$, we get, again in $\mathrm{Sp}_{U(1)}$,

$$\begin{aligned}\widetilde{\mathrm{Jac}}_m\colon \mathrm{MTSU}(m) &\xrightarrow{\eta\wedge\chi_m} \mathrm{TMF}\otimes B\mathrm{SU}(m)^{\overline{V}_m\boxtimes_{\mathbb{C}}\overline{L}+mL} \\ &\overset{\sigma}{\simeq} \mathrm{TMF}\otimes B\mathrm{SU}(m)_+\otimes S^{mL} \to \mathrm{TMF}\otimes S^{mL}\end{aligned}$$

by the $U(1)$-equivariant sigma orientation $\sigma$ and the collapse map $B\mathrm{SU}(m)_+ \to S^0$. The Jacobi orientation

$$\mathrm{Jac}_m\colon \mathrm{MTSU}(m) \to \mathrm{TJF}_m$$

is defined by passing to $U(1)$-fixed points.

**Theorem 3.6.** *The topological elliptic genera*

$$\bigoplus_{m\in\mathbb{Z}_{\geq 0}} \mathrm{Jac}_{U(1)_m}\colon \bigoplus_{m\in\mathbb{Z}_{\geq 0}} \mathrm{MTSU}(m) \to \bigoplus_{m\in\mathbb{Z}_{\geq 0}} \mathrm{TJF}_m$$

*assemble to an $E_1$-map. In particular, the stable topological elliptic genus*

$$\mathrm{Jac}\colon \mathrm{MSU} \to \mathrm{TJF}$$

*is $E_1$-multiplicative.*

*Proof.* In light of Lemma 3.4, it remains to show that the following map is $E_1$:

$$\sigma\colon \bigoplus_{m\in\mathbb{Z}_{\geq 0}} \left(\mathrm{TMF}\otimes BSU(m)^{\overline{V}_m\boxtimes_{\mathbb{C}}\overline{L}+mL}\right)^{U(1)} \simeq \bigoplus_{m\in\mathbb{Z}_{\geq 0}} \left(\mathrm{TMF}\otimes BSU(m)_+\otimes S^{mL}\right)^{U(1)}$$

This follows from the fact that both arrows in the following diagram are $E_\infty$.

$$BU(1)\times\left(\coprod_{m\in\mathbb{Z}_{\geq 0}} BSU(m)\right) \to BU(1)\times BSU \xrightarrow[\text{Def. 3.3}]{\mathfrak{s}} BU\langle 6\rangle.$$

Here, the first and second terms are $E_\infty$ algebras in the overcategory $\mathrm{Spaces}_{/BU(1)}$, and the last term is that in Spaces. The first arrow is the inclusion $BSU(m) \hookrightarrow BSU$ for the each component, so that the composition above gives the $BU\langle 6\rangle$-structure $\mathfrak{s}$ on $\overline{L}\boxtimes_{\mathbb{C}}\overline{V}_m$ for each $m$. The first map is obviously $E_1$. The second map is $E_\infty$, as verified in Definition 3.3. This completes the proof. ☐

The classical and topological elliptic genera are compatible in the sense that the following diagram commutes:

$$\begin{array}{ccc} \pi_*\,\mathrm{MTSU}(m) & \xrightarrow{\mathrm{Jac}_m} & \pi_*\mathrm{TJF}_m \\ \downarrow{\scriptstyle \mathrm{jac}_{\mathrm{clas},m}} & & \downarrow{\scriptstyle e_{\mathrm{JF}}} \\ \mathrm{jF}_m & \hookrightarrow & \mathrm{JF}_m \end{array} \tag{3.7}$$

3.4. **Two-variable elliptic genera à la Ando–French–Ganter.** In [AFG08], a variant of a stable topological elliptic genus is defined as a map

$$\sigma^{\#}\colon \mathrm{MSU} \to E^{\mathbb{CP}^{\infty}_{-\infty}},$$

naturally constructed for elliptic cohomology theories $E$. Here $\mathbb{CP}^{\infty}_{-\infty}$ is considered as the pro-spectrum $\{\mathbb{CP}^{\infty}_{-m} = (\mathbb{CP}^{\infty})^{-mL})\}_{m\geq 0}$, and the mapping spectrum in the target is by definition

$$E^{\mathbb{CP}^{\infty}_{-\infty}} = \operatorname*{colim}_m E^{\mathbb{CP}^{\infty}_{-m}}.$$

For $E = \mathrm{TMF}$, this mapping spectrum can be considered a Tate cohomology version of our TJF. Indeed, there is a cofiber sequence of pro-spectra

$$P_\infty \to D(\mathbb{CP}^{\infty}_{-\infty}) \to \Sigma^2 D(P_\infty)$$

where $P_\infty$ appeared in Theorem 2.4. Upon tensoring with TMF, gives a cofiber sequence

$$\mathrm{TJF} \to \mathrm{TMF}^{\mathbb{CP}^{\infty}_{-\infty}} \to \Sigma^2 D_{\mathrm{TMF}}\mathrm{TJF}.$$

We will start by describing the variant of $\sigma^{\#}$ we need. In [AFG08], a map

$$\mathrm{MSU} \to MU\langle 6\rangle^{\mathbb{CP}^{\infty}_{-\infty}}$$

is described, where $MU\langle 6\rangle$ is the bordism spectrum of manifolds with a lift of the classifying map of the stable normal bundle to the 5-connected cover $BU\langle 6\rangle$ of $BU$. We will refine this as follows:

$$\mathrm{MTSU}(m) \otimes \mathbb{CP}^0_{-m} \simeq B\mathrm{SU}(m)^{-\overline{V}_m} \otimes (\mathbb{CP}^m)^{-mL} = (B\mathrm{SU}(m) \times \mathbb{CP}^m)^{\overline{V}_m \boxtimes_{\mathbb{C}} \overline{L} - V_m \boxtimes_{\mathbb{C}} L}$$
$$\xrightarrow{\chi(V_m \boxtimes_{\mathbb{C}} L)} (B\mathrm{SU}(m) \times \mathbb{CP}^m)^{\overline{V}_m \boxtimes_{\mathbb{C}} \overline{L}} \xrightarrow{\mathfrak{s}} MU\langle 6\rangle$$

where we use the $BU\langle 6\rangle$-structure $\mathfrak{s}$ in Definition 3.3. The adjoint of the above map, together with the classical $\sigma$-orientation $\mathrm{MU}\langle 6\rangle \to \mathrm{MO}\langle 8\rangle \to \mathrm{tmf}$, gives the desired orientation

$$\widehat{\mathrm{jac}}_m\colon \mathrm{MTSU}(m) \to \mathrm{tmf} \otimes D(\mathbb{CP}^0_{-m}). \tag{3.8}$$

This is a destabilized version of the orientation $\sigma^{\#}$ of [AFG08] we will use in Section 4. Taking $m = \infty$, we get

$$\widehat{\mathrm{jac}}\colon \mathrm{MSU} = \mathrm{MTSU} \to \mathrm{tmf} \otimes D(\mathbb{CP}^0_{-\infty}),$$

which is an $E_\infty$ map by Definition 3.3. We also remark that the analogous $E_\infty$-ness of $\sigma^{\sharp}$ can also be verified in the same way and appears in [CL25].

## 4. Connective topological elliptic genera

In this section, we prove Thm.. 1.1, i.e. we construct a connective version of topological elliptic genera

$$\mathrm{jac}_m\colon \mathrm{MTSU}(m) \to \mathrm{tjF}_m.$$

which factor $\mathrm{Jac}_m\colon \mathrm{MTSU}(m) \to \mathrm{TJF}_m$ through the localization map $\mathrm{tjF}_m \to \mathrm{TJF}_m$.

As input for this construction, we need both versions of topological elliptic genera constructed in Section 3: the genus $\mathrm{Jac}_m\colon \mathrm{MTSU}(m) \to \mathrm{TJF}_m$ of [LY24] and the two-variable genus $\widehat{\mathrm{jac}}_m$ (3.8) essentially constructed in [AFG08].

Consider the following diagram, in which all the solid arrows are constructed:

(4.1)

$$\begin{array}{ccccc}
\mathrm{MTSU}(m) & & \xrightarrow{\widehat{\mathrm{jac}}_m} & & \\
{\scriptstyle \mathrm{jac}_m}\,\vdots & & & & \\
\mathrm{tjF}_m = \mathrm{tmf}\otimes P_m & \longrightarrow & \mathrm{tmf}\otimes D(\mathbb{CP}^0_{-m}) & \longrightarrow & \Sigma^2\mathrm{tmf} \\
\downarrow & & \downarrow & & \downarrow \\
\mathrm{TJF}_m \simeq \mathrm{TMF}\otimes P_m & \xrightarrow{\mathrm{id}_{\mathrm{TMF}}\otimes f_m} & \mathrm{TMF}\otimes D(\mathbb{CP}^0_{-m}) & \longrightarrow & \Sigma^2\mathrm{TMF}.
\end{array}$$

with a curved arrow $\mathrm{Jac}_m\colon \mathrm{MTSU}(m) \to \mathrm{TJF}_m$.

Here, the middle and bottom rows are fiber sequences induced from (2.5).

**Lemma 4.2.** *(The solid part of) diagram (4.1) commutes.*

*Proof.* The proof is analogous to [LY24, Section 4.3]. The map $\mathrm{id}_{\mathrm{TMF}}\otimes f_m\colon \mathrm{TMF}\otimes P_m \to \mathrm{TMF}\otimes D(\mathbb{CP}^0_{-m})$ is identified with

$$\begin{aligned}
\mathrm{TJF}_m &= F(S^{-mL}, \mathrm{TMF})^{U(1)} \\
&\to F(S((m+1)L)_+\otimes S^{-mL}, \mathrm{TMF})^{U(1)} \simeq \mathrm{TMF}\otimes D(\mathbb{CP}^0_{-m})
\end{aligned}$$

induced by the $U(1)$-equivariant collapse $S((m+1)L)_+ \to S^0$. The claim follows from the observation that the definition of $\mathrm{Jac}_m$ (Subsection 3.3) coincides with that of $\widehat{\mathrm{jac}}_m$ (Subsection 3.4) after inverting $\Delta^{-24}$, noting that $\mathbb{CP}^m_+ = S((m+1)L)_+/U(1)$. □

*Proof of Thm. 1.1.* We will show that there is a morphism

$$\mathrm{jac}\colon \mathrm{MTSU} = \mathrm{MSU} \to \mathrm{tjF}$$

which for $m=\infty$ makes Diagram (4.1) commute up to homotopy.

The general case (for finite $m \geq 0$) follows from this by considering the commutative diagram

$$\begin{array}{ccccc}
\mathrm{MTSU}(m) & \overset{\mathrm{jac}_m}{\dashrightarrow} & & \xrightarrow{\widehat{\mathrm{jac}}_m} & \\
 & & \mathrm{tjF}_m := \mathrm{tmf}\otimes P_m & \xrightarrow[\mathrm{id}_{\mathrm{tmf}}\otimes f_m]{} & \mathrm{tmf}\otimes D(\mathbb{CP}^0_{-m}) \\
\downarrow{\scriptstyle m\to\infty} & & \downarrow & & \downarrow \\
\mathrm{MTSU} & \xrightarrow{\mathrm{jac}} & \mathrm{tjF}\simeq \mathrm{tmf}\otimes P_\infty & \xrightarrow{\mathrm{id}_{\mathrm{tmf}}\otimes f_\infty} & \mathrm{tmf}\otimes D(\mathbb{CP}^0_{-\infty}),
\end{array}$$

where the solid part of the diagram commutes by the assumption on jac above. Moreover, the lower right square is a pullback diagram, so we get an induced dotted arrow $\mathrm{jac}_m$ as desired, and the compatibility across different values of $m$ comes for free.

To show the stable case, let us denote by $X$ the pullback

$$\begin{array}{ccc}
X & \longrightarrow & \mathrm{tmf}\otimes D(\mathbb{CP}^0_{-\infty}) \\
\downarrow & & \downarrow \\
\mathrm{TJF}\simeq \mathrm{TMF}\otimes P_\infty & \longrightarrow & \mathrm{TMF}\otimes D(\mathbb{CP}^0_{-\infty}).
\end{array}$$

By Lemma 4.2, the morphisms $\widehat{\mathrm{jac}}$ and Jac produces a morphism

$$\mathrm{jac}' \colon \mathrm{MTSU} \to X.$$

On the other hand, again by the universality of the pullback, we get a map $\phi\colon \mathrm{tjF} \to X$, and since the rows of (4.1) are fiber sequences, we get a fiber sequence

$$\mathrm{tjF} \xrightarrow{\phi} X \to \Sigma^1\mathrm{TMF}/\mathrm{tmf}.$$

The existence of the desired map jac is equivalent to the existence of a lift up to homotopy of $\mathrm{jac}'$ along the first map in the exact sequenc

$$[\mathrm{MTSU}, \mathrm{tjF}] \xrightarrow{\phi} [\mathrm{MTSU}, X] \to \left[\mathrm{MTSU}, \Sigma^1\mathrm{TMF}/\mathrm{tmf}\right].$$

We will show that

$$\left[\mathrm{MTSU}, \Sigma^1\mathrm{TMF}/\mathrm{tmf}\right] = 0, \tag{4.3}$$

from which the claim immediately follows. To show this, we use the mixed-case version Tmf in the map

$$\mathrm{TMF}/\mathrm{tmf} \to \mathrm{TMF}/\mathrm{Tmf} \simeq \mathrm{KO}((q))/\mathrm{KO}[[q]],$$

where the first map is the standard inclusion tmf $\to$ Tmf, and the latter equivalence follows from the fact [Lur09, Section 4.3] that the compactification $\overline{\mathcal{M}}$ of the elliptic moduli is formed by the pushout of $\mathcal{M}$ and $\operatorname{Spec} \mathbb{G}_m[[q]] /\!\!/ C_2$ along the Tate moduli $\operatorname{Spec} \mathbb{G}_m((q)) /\!\!/ C_2$. Since the map tmf $\to$ Tmf is $(-20)$-connective and $\mathrm{MTSU}(m)$ is connective, we get

$$\begin{aligned}\left[\mathrm{MTSU}, \Sigma^1\mathrm{TMF}/\mathrm{tmf}\right] &\simeq \left[\mathrm{MTSU}, \Sigma^1\mathrm{KO}[[q]]/\mathrm{KO}((q))\right] \\ &= \left[\mathrm{MTSU}, \Sigma^1\mathrm{KO}\right] \otimes (\mathbb{Z}[[q]]/\mathbb{Z}((q))).\end{aligned}$$

Furthermore, we have

$$\left[\mathrm{MTSU}, \Sigma^1\mathrm{KO}\right] \simeq \left[B\mathrm{SU}, \Sigma^1\mathrm{KO}\right]$$

by the Thom isomorphism, using the fact that the fundamental representation of $\mathrm{SU}(m)$ is KO-oriented.

The cofiber sequence $\Sigma\mathrm{KO} \xrightarrow{\eta} \mathrm{KO} \to \mathrm{KU}$ induces a long exact sequence

$$\cdots \to \mathrm{KU}^{-1}(B\mathrm{SU}) \to \mathrm{KO}^1(B\mathrm{SU}) \to \mathrm{KO}^0(B\mathrm{SU}) \to \cdots$$

Since $\eta \simeq 0$ in $B\mathrm{SU}$, the second map is zero. Since $B\mathrm{SU}$ only has cells in even dimensions, the Atiyah–Hirzebruch spectral sequence implies that also $\mathrm{KU}^{-1}(B\mathrm{SU}) = 0$ and thus $\mathrm{KO}^1(B\mathrm{SU}) = 0$. Thus (4.3) follows, completing the proof of Thm. 1.1. □

## 5. The ring $\pi_*\mathrm{tjF}$

In this section, we analyze the Adams-Novikov spectral sequence for tjF. Although we have not constructed a multiplication on tjF, we can find a natural multiplicative structure on the spectral sequence. The structure of the spectral sequence is determined in Corollary 5.7, and the resulting computation of $\pi_*\mathrm{tjF}$ as a ring is in Corollary 5.8.

We use the following conventions for modular forms. We denote by

$$\mathrm{mf} \cong \mathbb{Z}[c_4, c_6, \Delta]/(c_4^3 - c_6^2 - 1728\Delta)$$

the ring of integral, holomorphic modular forms (i.e. having integral Fourier coefficients in the variable $q = \exp(2\pi i\tau)$), and by

$$\mathrm{MF} := \mathrm{mf}[\Delta^{-1}]$$

the ring of integral, weakly homomorphic modular forms (i.e., holomorphic away from the cusp).

### 5.1. **Derived Jacobi forms.**

**Definition 5.1.** The bigraded ring of derived weakly holomorphic Jacobi forms is defined as

$$\mathrm{DJF}_{*,s} = H^s((\mathcal{E}^{\mathrm{or}})^\heartsuit - \{0\}; \pi_*(\mathcal{O}_{\mathcal{E}^{\mathrm{or}}-\{0\}}))$$

i.e. the $E^2$-terms of the descent spectral sequences, with $\mathrm{DJF}_{*,0}$ being identified with the ring of integral (even) Jacobi forms.

We wish to consider connective versions of the above, which for $s = 0$ coincide with weak Jacobi forms. As in [BM25], this can be accomplished with replacing the stack $\mathcal{M}_{\mathrm{ell}}^\heartsuit$ with the (classical) stack $\mathcal{M}_{\mathrm{Weier}}$ of cubic curves, which locally have a presentation by a Weierstrass equation. We denote the universal affine cubic over it by $\mathcal{W}$. We then have a pullback diagram

$$\begin{array}{ccc} (\mathcal{E}^{\mathrm{or}})^\heartsuit - \{0\} & \longrightarrow & \mathcal{W} \\ \downarrow & & \downarrow p \\ \mathcal{M}_{\mathrm{ell}}^\heartsuit & \longrightarrow & \mathcal{M}_{\mathrm{Weier}}, \end{array}$$

where the bottom map is the inclusion of the open substack $\Delta^{-1}\mathcal{M}_{\mathrm{Weier}}$.

The stack $\mathcal{M}_{\mathrm{Weier}}$ is not a Deligne-Mumford stack, nor is it the heart of a spectral stack for all we know. However, it has a flat affine cover by the spectrum of the ring

$$A = \mathbb{Z}[a_1, a_2, a_3, a_4, a_6],$$

as does the universal cubic $\mathcal{W}$ by

$$A_{\mathcal{W}} = A[x, y]/E(a_1, \ldots, a_6, x, y) \xrightarrow{\pi} \mathcal{W},$$

where $E$ is the Weierstrass equation (2.10)

$$E(a_1, \ldots, a_6, x, y) = y^2 + a_1xy + a_3y - (x^3 + a_2x^2 + a_4x + a_6).$$

Thus both $\mathcal{M}_{\mathrm{Weier}}$ and $\mathcal{W}$ are fppf stacks, whose cohomology with coefficients in the structure sheaf we denote by dmf and djF, respectively, the rings of derived modular forms and derived weak Jacobi forms, respectively. If we let

$$\mathrm{Spec}\,\Gamma_{\mathcal{W}} = \mathrm{Spec}\,A_{\mathcal{W}} \otimes_{\mathcal{W}} \mathrm{Spec}\,A_{\mathcal{W}}$$

then $\Gamma_{\mathcal{W}} = A[r, s, t]$, and $\mathcal{W}$ is the stack quotient of the Hopf algebroid $(A_{\mathcal{W}}, \Gamma_{\mathcal{W}})$.

**Definition 5.2.** The ring of weak derived Jacobi forms is defined as

$$\mathrm{djF}_{*,s} = H^s(\mathcal{W}; p^*\omega^*),$$

where $p\colon \mathcal{W} \to \mathcal{M}_{\mathrm{Weier}}$ is the canonical projection.

The computation of $\mathrm{djF}_{**}$ was carried out completely in [BM25]:

**Fact 5.3** ( [BM25, Theorem 4.5 and 4.7]). *We have*

$$\mathrm{djF}_{*,*} \cong \mathbb{Z}[b_2, b_3, b_4, b_8, h_1]/(2h_1, b_3h_1, 4b_8 + b_4^2 - b_2b_3^2)$$

*with*

$$|b_i| = (2i, 0);\ h_1 = (1, 1).$$

*where* $|b_i| = 2i$. *The* $(\mathrm{dmf}_{**})$*-module structure is given by*

| | | |
|---|---|---|
| $c_4 \mapsto b_2^2 - 24b_4$ | $c_6 \mapsto -b_2^3 + 36b_2b_4 - 216b_3^2$ | $\Delta \mapsto -b_2^2b_8 - 8b_4^3 - 27b_3^4 + 9b_2b_3^2b_4$ |
| $h_1 \mapsto h_1$ | $\alpha_1 \mapsto 0$ | $\delta \mapsto b_2h_1$ |
| $\epsilon \mapsto 0$ | $B_1 \mapsto 0$ | $\kappa \mapsto 0$ |
| $\bar{\kappa} \mapsto 0.$ | | |

We recover the structure of jF as the subring of elements of degree $(*, 0)$.

5.2. **Comparison with Gritsenko's computation and identification of the generators.** Unstably, by [Gri20, Theorem 2.7], we have the generator–relation expression,

$$\mathrm{jF}_\bullet = \mathrm{mf}[\phi_{-1,\frac{1}{2}}, \phi_{0,1}, \phi_{0,\frac{3}{2}}, \phi_{0,2}, \phi_{0,4}, E_{4,1}, E_{4,2}, E_{4,3}, E_{6,1}, E_{6,2}, E'_{6,3}]/\sim,$$

where for the relation $\sim$ we refer to loc. cit.

The notation $f_{k,m}$ denotes an element of weight $k$ and index $m$, so that $f_{k,m} \in \mathrm{jF}_{2m}$ with degree $2k + 4m$. (Note that we follow the notation in the literature of not doubling the (half-integer) index for these classes, as we otherwise do in this article.)

The generators in Gritsenko's description each have simple expressions in terms of theta, eta, and Weierstrass $\wp$-functions together with functions

$$\xi_{ab}(z, q) = \theta_{ab}(z, q)/\theta_{ab}(0, q).$$

They stabilize to generators in our description, Fact 5.3, as summarized in Table 1.

| Generator in [Gri20] | analytic expression | stabilization in jF |
|---|---|---|
| $\phi_{-1,\frac{1}{2}}$ | $\frac{\theta_{11}(z,\tau)}{\eta^3(\tau)}$ | $a$ |
| $\phi_{0,1}$ | $-\frac{3}{\pi^2}\wp(z,\tau)\frac{\theta_{11}^2(z,\tau)}{\eta^6(\tau)}$ | $b_2$ |
| $\phi_{0,\frac{3}{2}}$ | $\frac{\theta_{11}(2z,\tau)}{\theta_{11}(z,\tau)}$ | $b_3$ |
| $\phi_{0,2}$ | $2(\xi_{00}^2 + \xi_{01}^2 + \xi_{10}^2)(2z, \tau)$ | $b_4$ |
| $\phi_{0,4}$ | $\frac{\theta_{11}(3z,\tau)}{\theta_{11}(z,\tau)}$ | $b_8$ |

TABLE 1. Gritsenko's generators and their stabilization

5.3. **Descent and Adams-Novikov spectral sequences.** The aim of this section is to construct a descent spectral sequence for connective tjF and show that it is isomorphic with the MU-based Adams–Novikov spectral sequence. The arguments are completely analogous to the corresponding results about tmf, cf. [DFHH14, Chapter 9.3].

**Definition 5.4** (Ravenel's filtration of MU)**.** The spectrum $X(n)$ is the Thom spectrum of the vector bundle on $\Omega\mathrm{SU}(n)$ classified by

$$\Omega\mathrm{SU}(n) \hookrightarrow \Omega\mathrm{SU} \xrightarrow{\simeq} BU,$$

where the last map is the Bott periodicity equivalence.

The crucial feature of these spectra is that for complex oriented spectra $E$, homotopy classes of homotopy multiplicative maps $X(n) \to E$ are in one-to-one correspondence with coordinates modulo degree $n+1$ on the formal group $\mathbb{G}_E$ associated to $E$ [DFHH14, Ch. 9]. As a consequence, since an elliptic curve together with a coordinate modulo degree 5 exactly fixes a Weierstrass parametrization, we have

$$\pi_*(\mathrm{tmf} \otimes X(4)) = A = \mathbb{Z}[a_1, a_2, a_3, a_4, a_6],$$
$$\pi_*(\mathrm{tmf} \otimes X(4) \otimes X(4)) = \Gamma = A[r, s, t],$$

and

$$\pi_*(\mathrm{tmf} \otimes X(4)^{\otimes n+1}) \cong \Gamma \otimes_A \cdots \otimes_A \Gamma \quad (n \text{ factors.})$$

Thus the descent spectral sequence is exactly the $X(4)$-based Adams-Novikov spectral sequence for the spectrum TMF, regardless of the fact that $X(4)$ is not actually complex oriented, and it is reasonable to also call the tmf-based version a descent spectral sequence. Concretely, it is the Bousfield–Kan spectral sequence associated with the cosimplicial resolution

$$\mathrm{tmf} \to \mathrm{tmf} \otimes X(4) \rightrightarrows \mathrm{tmf} \otimes X(4) \otimes X(4) \cdots .$$

The inclusion map $X(4) \to X(\infty) = MU$ gives a map from the resolution above to the MU-based Adams–Novikov resolution and thus a map of spectral sequences from the descent spectral sequence to the Adams–Novikov spectral sequence for tmf. In fact, it induces an isomorphism from the $E_2$-term on since

$$\pi_*(\mathrm{MU} \otimes \mathrm{tmf}) \cong \pi_*(X(4) \otimes \mathrm{tmf})[x_5, x_6, \dots].$$

Just like $\operatorname{Spec} \pi_*(X \otimes \mathrm{tmf})$, also $\operatorname{Spec} \pi_*(MU \otimes \mathrm{tmf}) \to \mathcal{M}_{\mathrm{Weier}}$ is a flat cover, albeit a less efficient one, so that

$$\begin{aligned}\mathrm{dmf}_{**} = \mathrm{Ext}^{**}_{(A,\Gamma)}(A, A) &\cong \mathrm{Ext}^{**}_{(X(4)_*, X(4)_* X(4))}(X(4)_*, X(4)_* \mathrm{tmf}) \\ &\cong \mathrm{Ext}^{**}_{(\mathrm{MU}_*, \mathrm{MU}_* \mathrm{MU})}(\mathrm{MU}_*, \mathrm{MU}_* \mathrm{tmf}).\end{aligned}$$

For topological Jacobi forms, we have a pullback diagram of stacks

$$\begin{array}{ccc} \operatorname{Spec} X(4)_* \mathrm{tjF} & \longrightarrow & \mathcal{W} \\ \downarrow & & \downarrow \\ \operatorname{Spec} X(4)_* \mathrm{tmf} & \longrightarrow & \mathcal{M}_{\mathrm{Weier}}, \end{array}$$

and since the lower horizontal map is a flat cover, so is the upper one, and thus so is the composite $\operatorname{Spec} \mathrm{MU}_* \mathrm{tjF} \to \operatorname{Spec} X(4)_* \mathrm{tjF} \to \mathcal{W}$. The homology spectral sequence associated to this cover can reasonably be called the descent spectral sequence for connective tjF, isomorphic to the periodic descent spectral sequence (2.9) upon inverting $\Delta$. Since $\operatorname{Spec} \mathrm{MU}_* \mathrm{tjF} \to \mathcal{W}$ is a flat cover, it is isomorphic to the Adams–Novikov spectral sequence. We obtain:

**Lemma 5.5.** *The Adams–Novikov spectral sequence for* tjF *maps to the descent spectral sequence* (2.9) *for* TJF*, and this map becomes an isomorphism upon inverting $\Delta$. Moreover, the Adams–Novikov $E_2$-term for* tjF *is isomorphic to* djF*:*

$$E_2^{\mathrm{ANSS}}(\mathrm{tjF}) \simeq \mathrm{djF}_{*,*} .$$

The structure of descent (or Adams–Novikov) spectral sequence for tjF was completely determined in [BM25]. We have **Fact 5.6** ( [BM25, Section 6]). *The Adams–Novikov spectral sequence for* tjF *collapses at $E_4$. It has no nontrivial $d^2$-differentials, and all $d^3$-differentials are generated by the Leibniz rule and*

$$d^3(b_2) = h_1^3.$$

By Lemma 5.5 and Fact 5.6, we get

**Corollary 5.7.** *The localization map $E^{\mathrm{ANSS}}(\mathrm{tjF}) \to E^{\mathrm{ANSS}}(\mathrm{TJF})$ exibits $E^{\mathrm{ANSS}}(\mathrm{tjF})$ as a multiplicative sub-spectral sequence of $E^{\mathrm{ANSS}}(\mathrm{TJF})$; namely, we have $E_2^{\mathrm{ANSS}}(\mathrm{tjF}) = \mathrm{djF}_{*,*}$, and all $d^3$-differentials are generated by Leibniz rule and $d^3(b_2) = h_1^3$, and we have $d^r = 0$ for $r \geq 4$.*

Note that this is a nontrivial result based on computations, given that we are not providing tjF with a multiplicative structure. In particular, it provides $\pi_*\mathrm{tjF}$ with a ring structure.

5.4. **Determination of the ring $\pi_*\mathrm{tjF}$.** By the computations in [BM25], we get

**Corollary 5.8** (of [BM25, Corollary 6.1]). *The map $\pi_*\mathrm{tjF} \to \pi_*\mathrm{TJF}$ exhibits $\pi_*\mathrm{tjF}$ as a subring of $\pi_*\mathrm{TJF}$. The ring of weak topological Jacobi forms, $\pi_*\mathrm{tjF}$, is generated by the following classes:*

- *Classes $x_i$ of degree $2i$ for $i = 2, 3, 4, 5, 6, 8$;*
- *A class $x_4'$ of degree 8,*
- *A class $\eta$ of degree 1.*

*The relations between the generators $x_i, x_4'$ are given as the kernel of edge map*

$$\mathbb{Z}[x_i, x_4'] \to \mathbb{Z}[b_2, b_3, b_4, b_8]/(4b_8 + b_4^2 - b_2 b_3^2) = \mathrm{jF}$$

*mapping*

$$x_2 \mapsto 2b_2 \quad x_3 \mapsto b_3 \quad x_4 \mapsto b_2^2 \quad x_4' \mapsto b_4 \quad x_5 \mapsto b_2 b_3 \quad x_6 \mapsto b_2 b_4 \quad x_8 \mapsto b_8$$

*In addition, the relations*

$$2\eta = 0, \quad \eta^3 = 0, \quad x_i\eta = 0 \text{ for all } i \neq 4, 8, \quad x_4'\eta = 0$$

*hold, and these complete the relations among the generators of $\pi_*\mathrm{tjF}$.*

The edge homomorphisms give us the commutative diagram of graded rings

$$\begin{array}{ccc} \pi_*\mathrm{tjF} & \xrightarrow{e_{\mathrm{jF}}} & \mathrm{jF} \\ \downarrow & & \downarrow \\ \pi_*\mathrm{TJF} & \xrightarrow{e_{\mathrm{JF}}} & \mathrm{JF}. \end{array}$$

## 6. The homotopy of MSU

By the theorems of Quillen and Lazard, the coefficients of MU are isomorphic to the Lazard ring

$$\mathrm{MU}_* \cong \mathbb{Z}[\alpha_1, \alpha_2, \dots]; \quad |\alpha_i| = 2i,$$

over which the universal formal group law is defined. The polynomial generators $\alpha_i$ can be taken as any set of stably complex $2i$-dimensional manifolds $M_i$ as long as their Milnor numbers $s_i(M_i) \in H^{2i}(M_i; \mathbb{Z}) \cong \mathbb{Z}$ satisfy

$$s_i(M_i) = \pm m_i,$$

where

$$m_i = \begin{cases} p; & i+1 \text{ is a prime power;} \\ 1; & \text{otherwise.} \end{cases} \tag{6.1}$$

The Milnor number $s_n$ is the characteristic number associated to the power sum $x_1^n + \cdots + x_n^n$ of the Chern roots of the tangent bundle.

Unlike MU, the spectrum MSU of SU-bordism has interesting 2-torsion and has been studied since the 1960s. The graded ring structure of $\mathrm{MSU}_*$ is completely known [Sto15, Chapter X]. Our strategy to show the surjectivity of the topological elliptic genus in homotopy, Thm. 1.2, is to show surjectivity on the Adams–Novikov $E_2$-terms of MSU and tjF, and trace the map through the spectral sequence. Thus we not only need to know $\mathrm{MSU}_*$, but also the structure of the Adams–Novikov spectral sequence converging to it.

6.1. **Review on the structure of** $\mathrm{MSU}_*$**.** For details of the content of this subsection, we refer the reader to [Sto15, Chapter X] or [CLP19].

To understand MSU, consider the fibration

$$S^3 \xrightarrow{\eta} \mathbb{CP}^1 = S^2 \xrightarrow{\iota} BU(1) = \mathbb{CP}^\infty,$$

where the map $\iota$ is the standard inclusion and $\eta$ is the Hopf map. Let $W$ be the pullback in

$$\begin{array}{ccc} W & \longrightarrow & BU \\ \big\downarrow & & \big\downarrow{\scriptstyle \det} \\ S^2 & \xrightarrow{\iota} & BU(1) \end{array}$$

and $MW$ the Thom spectrum associated to $W \to BU$. Geometrically, this is the bordism theory of stably almost complex manifolds with a lift of the classifying map of their determinant line bundles to $\mathbb{CP}^1$. The bordism theory $MW$ sits in between MSU and MU and is useful to compute $\mathrm{MSU}_*$. Indeed, Stong [Sto15, Chapter X, from page 262] constructs a split short exact sequence

$$0 \to MW_* \to \mathrm{MU}_* \xrightarrow{d} \mathrm{MU}_{*-4} \to 0,$$

where the map $d \in \mathrm{MU}^4\mathrm{MU}$ is the image under the Thom isomorphism from $\mathrm{MU}^4(BU) \cong \mathrm{MU}_*[\![c_1, c_2, \cdots]\!]$ of the element $-c_1^2(\det(V))$, where $V$ is the universal virtual bundle over $BU$.

Moreover, we have

$$MW \simeq \mathrm{MSU}/\eta,$$

the cone of multiplication with the Hopf map $\eta \in \pi_1 S$. Thus there is a fiber sequence

$$\Sigma^1 \mathrm{MSU} \xrightarrow{\eta} \mathrm{MSU} \xrightarrow{q} MW \xrightarrow{\partial} \Sigma^2 \mathrm{MSU}. \tag{6.2}$$

A slight complication is that although $\mathrm{MSU} \to \mathrm{MU}$ is a map of commutative ring spectra, $MW_*$ is not a subring of $\mathrm{MU}_*$. However, there is a twisted multiplication on $MW_*$ in such a way that $\mathrm{MSU}_* \to MW_*$ are ring maps. In this structure,

$$MW_* \cong \mathbb{Z}[\rho_1, \rho_3, \rho_4, \rho_5, \cdots]$$

is a polynomial algebra on generators $|\rho_i|$ in all degrees $2i$ except for 4. By [CLP19, Thm. 6.10], the generators can be chosen arbitrarily as long as they satisfy

$$s_i(\rho_i) = \pm m_i m_{i-1},$$

with $m_i$ as in (6.1).

A choice can be made such that

$$\rho_1 = [\mathbb{CP}^1], \; q \circ \partial(\rho_{2k}) = \rho_{2k-1} \text{ for } k \geq 2,$$

with $\delta\colon \mathrm{M}W \to \Sigma^2 \mathrm{MSU}$ and $q\colon \Sigma^2 \mathrm{MSU} \to \Sigma^2 \mathrm{M}W$ as in (6.2). We will fix such a (non-unique) choice of generators.

The Bockstein spectral sequence

$$E_1 = \mathrm{M}W_*[h_1] \Longrightarrow \mathrm{MSU}_*$$

has $d^1(\rho_1) = 2h_1$, $d^1(\rho_{2n}) = \rho_{2n-1}h_1$, and $d^3(\rho_1^2) = h_1^3$. The differential $d^1$ is not a derivation; instead, it satisfies $d^1(xy) = d^1(x)y \pm x d^1(y) \pm \rho_1 \frac{1}{h_1} d^1(x) d^1(y)$.

In fact, it is isomorphic to the Adams–Novikov spectral sequence, as we will now show. The Bockstein spectral sequence is associated to the resolution

$$\mathrm{M}W \xrightarrow{q\circ\partial} \Sigma^2\mathrm{M}W \xrightarrow{q\circ\partial} \Sigma^4\mathrm{M}W \to \cdots$$

and there is a map of resolutions

$$\begin{array}{ccccccccc}
\mathrm{MSU} & \longrightarrow & \mathrm{M}W & \xrightarrow{q\circ\partial} & \Sigma^2\mathrm{M}W & \xrightarrow{q\circ\partial} & \Sigma^4\mathrm{M}W & \longrightarrow & \cdots \\
\| & & \downarrow & & \downarrow{\scriptstyle \binom{\mathrm{id}}{0}} & & \downarrow{\scriptstyle \binom{\mathrm{id}}{0}} & & \\
\mathrm{MSU} & \longrightarrow & \mathrm{MU} & \xrightarrow{\binom{\partial}{d}} & \Sigma^2\mathrm{MU} \vee \Sigma^4\mathrm{MU} & \longrightarrow & \Sigma^4\mathrm{MU} \vee \Sigma^6\mathrm{MU} & \longrightarrow & \cdots
\end{array}$$

This produces a map from the Bockstein spectral sequence to the Adams–Novikov spectral sequence since the bottom row is an Adams–Novikov resolution. Moreover, it induces an isomorphism on the $E_2$-term as shown in [CLP19, Prop. 5.1].

**Corollary 6.3.** *The Adams–Novikov spectral sequence converging to* $\mathrm{MSU}_*$ *has*

$$E_2^{**} = \frac{\mathbb{Z}[h_1, B_2, B_3, B_4, \dots, C_8, C_{12}, C_{16}, \dots]}{(2h_1, B_{2n+1}h_1, B_{2n}^2 + 4C_{4n} - B_2 B_{2n-1})}$$

*with* $|h_1| = (1,2)$, $|B_n| = (0, 2n)$, $|C_{4n}| = (0, 8n)$ .

*In the* $E_1$*-term of the Bockstein spectral sequence, these classes are represented as follows:*

$$B_2 = \rho_1^2; \quad B_{2n+1} = \rho_{2n+1}; \quad B_{2n+2} = 2\rho_{2n+2} - \rho_1\rho_{2n+1}; \quad C_{4n} = \rho_1\rho_{2n-1}\rho_{2n} - \rho_{2n}^2.$$

*The spectral sequence collapses at* $E_4$. *The* $d^2$*-differentials are trivial, and all* $d^3$*-differentials follow from the Leibniz rule and*

$$d^3(B_2) = h_1^3.$$

We observe that there is no room for multiplicative extensions and conclude that the graded ring $\mathrm{MSU}_*$ is isomorphic to the cohomology ring of the differential graded algebra $(E_2^{**}, d^3)$ above (cf. [Sto15]). In the range $n \leq 16$, the additive generators of $\mathrm{MSU}_n$ in terms of the generators in Cor. 6.3 are listed in Table 2. Here, for an infinite cycle $x \in E_r$ in the Adams–Novikov spectral sequence, we denote by $[x] \in \pi_* R$ any choice of corresponding class represented by $x$. The multiplicative structure is read off from that. We use the conventional notation $\eta := [h_1]$ since it is the image of the generator $\eta \in \pi_1 S$.

| $n$ | $\pi_n\,\mathrm{MSU}$ | additive generators of $\pi_n\,\mathrm{MSU}/\mathrm{tors}$ | $\mathrm{tors}(\pi_n\,\mathrm{MSU})$ |
|---|---|---|---|
| 0 | $\mathbb{Z}$ | 1 | |
| 1 | $\mathbb{Z}/2$ | | $\eta$ |
| 2 | $\mathbb{Z}/2$ | | $\eta^2$ |
| 3 | 0 | | |
| 4 | $\mathbb{Z}$ | $[2B_2]$ | |
| 5 | 0 | | |
| 6 | $\mathbb{Z}$ | $[B_3]$ | |
| 7 | 0 | | |
| 8 | $\mathbb{Z}^2$ | $[B_2^2], [B_4]$ | |
| 9 | $\mathbb{Z}/2$ | | $\eta[B_2^2]$ |
| 10 | $\mathbb{Z}^2 \oplus \mathbb{Z}/2$ | $[B_2B_3], [B_5]$ | $\eta^2[B_2^2]$ |
| 11 | 0 | | |
| 12 | $\mathbb{Z}^4$ | $[2B_2^3], [B_3^2], [B_2B_4], [B_6]$ | |
| 13 | 0 | | |
| 14 | $\mathbb{Z}^4$ | $[B_2^2B_3], [B_2B_5], [B_4B_3], [B_7]$ | |
| 15 | 0 | | |
| 16 | $\mathbb{Z}^7$ | $[B_2^4], [B_2^2B_4], [B_2B_3^2], [B_2B_6], [C_8], [B_3B_5], [B_8]$ | |

TABLE 2. $\pi_n\,\mathrm{MSU}$ for $n \le 16$

**Proposition 6.4** ( [CLP19, Thm. 7.1]). *The following indecomposable elements have minimal positive Milnor numbers in* $\mathrm{MSU}_{2n}$*:*

- *for* $n = 2$, $[2B_2] = [K3]$ *with* $s_2([2B_2]) = -48$ *(K3 being a K3-surface);*
- *for* $n = 2i \ge 4$, $[B_{2i}]$ *with* $s_{2i}([B_{2i}]) = 2m_{2i}m_{2i-1}$ *with* $m_i$ *as in* (6.1)*;*
- *for* $n = 2i + 1$, $[B_{2i+1}]$ *with* $s_{2i+1}([B_{2i+1}]) = m_{2i+1}m_{2i}$*.*

6.2. **Elliptic genera of some generators of** $\mathrm{MSU}_*$**.** By (3.7) and (4.1), we have

$$\mathrm{jac}_{\mathrm{clas}} = e_{\mathrm{jF}} \circ \mathrm{jac}\colon\ \mathrm{MSU}_* \to \pi_*\mathrm{TJF} \to \mathrm{JF} \tag{6.5}$$

We can compute $\mathrm{jac}_{\mathrm{clas}}$ in terms of characteristic numbers. In low degrees, we get:

**Proposition 6.6.** *The elliptic genera restricted to* $\mathrm{MSU}_n$ *for* $n = 4, 6, 8$ *are given by*

$$\begin{aligned}
\mathrm{jac}_{\mathrm{clas}}|_{\mathrm{MSU}_4} &= \frac{c_2}{12} \cdot b_2 = -\frac{s_2}{24} \cdot b_2,\\
\mathrm{jac}_{\mathrm{clas}}|_{\mathrm{MSU}_6} &= \frac{c_3}{2} \cdot b_3 = \frac{s_3}{6} \cdot b_3,\\
\mathrm{jac}_{\mathrm{clas}}|_{\mathrm{MSU}_8} &= -\frac{s_4}{20} \cdot b_4 + \left(\frac{c_2^2}{240} - \frac{c_4}{720}\right) \cdot b_2^2.
\end{aligned}$$

*Proof.* The first two equations follow from the fact that for any closed SU-manifold $M$, we have

$$\mathrm{jac}_{\mathrm{clas}}(M) = \chi(M) + O(z)((q)) \quad (\chi \text{ being the Euler characteristic})$$

and

$$b_2 = 12 + O(z)((q)), \quad b_3 = 2 + O(z)((q)).$$

For the elliptic genus in degree 8, we use that it can be written in terms of Chern numbers as[4]

$$\mathrm{jac}_{\mathrm{clas}}|_{\mathrm{MSU}_8} = \left(\frac{1}{240}c_2^2 - \frac{1}{720}c_4\right) y^{\pm 2} + \left(-\frac{1}{60}c_2^2 + \frac{31}{180}c_4\right) y^{\pm 1} + \left(\frac{1}{40}c_2^2 + \frac{79}{120}c_4\right) + O(q).$$

On the other hand, we have

$$s_4 = 2c_2^2 - 4c_4. \quad \text{on } \mathrm{MSU}_*$$

Furthermore, we use (see, e.g. [Gri99, Section 1] together with Table 1 for the formulas)

$$b_2^2 = y^{\pm 2} + 20y^{\pm 1} + 102 + O(q)$$
$$b_4 = y^{\pm 1} + 4 + O(q)$$

Combining these completes the proof. □

**Corollary 6.7.** *The elliptic genera of the additive generators of* $\mathrm{MSU}_n$ *for* $n = 4, 6, 8$ *in Table 2 are given as follows.*

*(1)* $\mathrm{jac}_{\mathrm{clas}}([2B_2]) = 2b_2$*;*

*(2)* $\mathrm{jac}_{\mathrm{clas}}([B_3]) = b_3$*;*

*(3)* $\mathrm{jac}_{\mathrm{clas}}([B_2^2]) = b_2^2$*;*

*(4)* $\mathrm{jac}_{\mathrm{clas}}([B_4]) = -b_4 + 2Nb_2^2$ *for some* $N \in \mathbb{Z}$*.*

*Proof.* (1) and (2) follow by applying the formulas in Proposition 6.6 to those in Proposition 6.4. (3) follows from (1) and the multiplicativity of the elliptic genera. We similarly get from those two propositions that

$$\mathrm{jac}_{\mathrm{clas}}([B_4]) = -b_4 + nb_2^2 \quad \text{for some } n \in \mathbb{Z}.$$

To see that $n$ is even, consider $[B_2B_4] \in \pi_{12}\,\mathrm{MSU}$ (cf. Table 2). This element satisfies

$$\mathrm{jac}_{\mathrm{clas}}([B_2B_4]) = \frac{1}{2}\mathrm{jac}_{\mathrm{clas}}([2B_2]) \cdot \mathrm{jac}_{\mathrm{clas}}([B_4]) = -b_2b_4 - nb_2^3.$$

By the factorization (6.5), the above Jacobi form msut be contained in the image of $e_{\mathrm{JF}}$. By Fact 5.6, the image has index 2 and we conclude that $n$ is even. This verifies (4) and completes the proof. □

*Remark* 6.8. Corollary 6.7 implies in particular the surjectivity of

$$\mathrm{jac}_{\mathrm{clas}} \otimes \mathbb{Z}[\tfrac{1}{2}]\colon\ \mathrm{MSU}_*[\tfrac{1}{2}] \to \mathrm{jF} \otimes \mathbb{Z}[\tfrac{1}{2}] = \mathbb{Z}[\tfrac{1}{2}, b_2, b_3, b_4].$$

⌟

[4]For this computation, the authors used the sage package for calculation of elliptic genera developed by Kenta Kobayashi [Kob23], available in https://github.com/topostaro/EllipticGenus.

## 7. The surjectivity of the topological elliptic genus

The goal of this section is to show Thm. 1.2 by studying the map induced on Adams–Novikov spectral sequences by $\mathrm{jac}\colon \mathrm{MSU} \to \mathrm{tjF}$.

By Remark 6.8 and the fact that $\pi_*\mathrm{tjF}[\frac{1}{2}] \simeq \mathrm{jF}[\frac{1}{2}]$, it is enough to show the result 2-locally. Thus, in the section, all rings and spectra are implicitly localized at 2.

For a spectrum $R$ we denote by $E_r^{**}(R)$ the $E_r$-page of the Adams–Novikov spectral sequence for $R$.

By Corollary 5.7, the map

$$E_*(\mathrm{jac})\colon\ (E_r^{p,q}(\mathrm{MSU}), d^r) \to (E_r^{p,q}(\mathrm{tjF}), d^r)\,. \tag{7.1}$$

is multiplicative because the composition $\mathrm{MSU} \to \mathrm{tjF} \to \mathrm{TJF}$ is multiplicative, and that the map $\mathrm{tjF} \to \mathrm{TJF}$ induces an injection of spectral sequences.

The structures of the Adams–Novikov spectral sequences for MSU and tjF were determined in Corollaries 6.3 and 5.7. In particular, the $E_2$ pages are given by

$$E_2(\mathrm{MSU}) = \frac{\mathbb{Z}[h_1, B_2, B_3, B_4, \ldots, C_8, C_{12}, C_{16}, \ldots]}{(2h_1, B_{2n+1}h_1, B_{2n}^2 + 4C_{4n} - B_2B_{2n-1})}$$

and

$$E_2(\mathrm{tjF}) \simeq \mathrm{djF}_{*,*} \cong \frac{\mathbb{Z}[b_2, b_3, b_4, b_8, h_1]}{(2h_1, b_3h_1, 4b_8 + b_4^2 - b_2b_3^2)}.$$

As an immediate consequence of Corollary 6.3, we observe:

**Lemma 7.2.** *The generators $h_1, B_2, B_3, B_4, C_8$ of $E_2(\mathrm{MSU})$ generate a multiplicative sub-spectral sequence $(\widetilde{E}_*(\mathrm{MSU}), \widetilde{d}^*)$ of $(E_*(\mathrm{MSU}), d^*)$ whose $E_2$ page is*

$$\widetilde{E}_2(\mathrm{MSU}) = \mathbb{Z}[h_1, B_2, B_3, B_4, C_8]/(2h_1, B_3h_1, 4C_8 + B_4^2 - B_2B_3)$$

*with all nontrivial differentials generated under the Leibniz rule by $d^3(B_2) = h_1^3$.* □

**Proposition 7.3.** *The map (7.1) maps the sub-spectral sequence $\widetilde{E}(\mathrm{MSU})$ isomorphically to $E(\mathrm{tjF})$.*

*Proof.* Consider the map of $E_2$-pages,

$$E_2(\mathrm{jac})|_{\widetilde{E}(\mathrm{MSU})}\colon \widetilde{E}_2(\mathrm{MSU}) \to E_2(\mathrm{tjF}) \simeq \mathrm{djF}_{*,*}\,.$$

We already know that $h_1 \in \widetilde{E}_2(\mathrm{MSU})$ maps to $h_1 \in E_2(\mathrm{tjF})$ since they are coming from $E_2(S)$. We check the rest of the generators, which are all in $\widetilde{E}_2^{0,*}$. We use the commutative diagram

$$\begin{array}{ccc}
\mathrm{MSU}_* & \xrightarrow{\ \mathrm{jac}\ } & \pi_*\mathrm{tjF} \\
\downarrow{\scriptstyle e_{\mathrm{MSU}}} & \searrow^{\mathrm{jac}_{\mathrm{clas}}} & \downarrow{\scriptstyle e_{\mathrm{jF}}} \\
E_2^{0,*}(\mathrm{MSU}) & \xrightarrow{\ E_2(\mathrm{jac})\ } & \mathrm{djF}_{*,0} = \mathrm{jF}
\end{array}$$

where the vertical arrows are the edge homomorphisms of the spectral sequences. We use the fact that the diagonal arrow is the classical elliptic genus, which we already computed for some generators of $\mathrm{MSU}_*$ in Corollary 6.7. By that collorary and Corollary 6.3, we see that $E_2(\mathrm{jac})$ maps

$$\begin{aligned} B_2 &\mapsto b_2 \\ B_3 &\mapsto b_3, \\ B_4 &\mapsto -b_4 + 2Nb_2^2, \end{aligned}$$

$$C_8 = \frac{1}{4}\left(B_4^2 - B_2B_3^2\right) \mapsto \frac{1}{4}\left((-b_4+2Nb_2^2)^2 - b_2b_3^2\right) = -b_8 - Nb_2^2b_4 + N^2b_2^4.$$

Thus we obtain an isomorphism of $E_2$-pages and hence between all of the spectral sequences. □

*Proof of Thm. 1.2.* Combine Lemma 7.2 and Proposition 7.3. □

In particular, we can identify the image of the classical elliptic genera as in Corollary 1.3.

## 8. The invariance of the $U(1)$-topological elliptic genera under SU-flops

Totaro [Tot00] identifies the kernel of the classical elliptic genera with the ideal of $\pi_* \mathrm{MSU} \otimes \mathbb{Z}[\frac{1}{2}]$ generated by SU-*flops*. We expect an integral version of this statement to hold for our topological elliptic genera. As a first step, we show that topological elliptic genera are invariant under SU-flops.

**8.1. SU-flops.** We give a homotopy-theoretic formulation of Totaro's notion of "SU-flops" [Tot00].

**Definition 8.1** (twisted projective space $\widetilde{\mathbb{CP}}_{2,2}$)**.** We define a real 6-dimensional stable tangential $SU(4)$-manifold $\widetilde{\mathbb{CP}}_{2,2}$ with $U(2)\times U(2)$-action as follows. Given two 2-dimensional complex vector spaces $A$ and $B$, and fix identifications $A \simeq \mathbb{C}^2$ and $B \simeq \mathbb{C}^2$. Let

$$\widetilde{\mathbb{CP}}_{2,2} = \mathbb{CP}(A \oplus B^*) \cong \mathbb{CP}^3.$$

As explained in [Tot00, p.773], we have an isomorphism of *real* vector bundles over $\widetilde{\mathbb{CP}}_{2,2}$,

$$T(\widetilde{\mathbb{CP}}_{2,2}) \oplus \underline{\mathbb{R}}^2 \simeq A \otimes_{\mathbb{C}} \mathcal{O}(1) \oplus B \otimes_{\mathbb{C}} \mathcal{O}(-1), \tag{8.2}$$

where $\mathcal{O}(\pm 1)$ are the tautological complex line bundle and its dual on $\widetilde{\mathbb{CP}}_{2,2}$, and $\underline{\mathbb{R}}$ denotes the trivial line bundle. The right hand side of (8.2) is a SU(4)-vector bundle, giving $\widetilde{\mathbb{CP}}_{2,2}$ a tangential SU(4)-structure. We call $\widetilde{\mathbb{CP}}_{2,2}$ with this structure the *twisted complex projective space.*

The natural $U(A)\times U(B)$-action on $\widetilde{\mathbb{CP}}_{2,2}$ does not preserve this SU(4)-structure on the nose, but it does up to a twist by $-(A\oplus B)$.[5] By the identification $A \simeq \mathbb{C}^2$ and $B \simeq \mathbb{C}^2$ fixed in the beginning, we regard $\widetilde{\mathbb{CP}}_{2,2}$ as a stably tangential $SU(4)$-manifold with this twisted $U(2)\times U(2)$-action.

The corresponding bordism class in the $U(2)\times U(2)$-equivariant Thom spectrum MTSU(4) is

$$[\widetilde{\mathbb{CP}}_{2,2}] \in \pi_6\, \mathrm{MTSU}(4)\left[-\left(\overline{V}_{U(2)} \oplus \overline{V}_{U(2)}\right)\right]^{U(2)\times U(2)}. \tag{8.3}$$

[5] The sign convention is made so that (8.3) holds.

**Definition 8.4** ($SU$-bordism modulo flops)**.** For each $m \geq 4$, define $\mathrm{flop}_m$ to be the following composition:

$$\begin{aligned}\mathrm{flop}_m \colon \Sigma^6 \, \mathrm{MTSU}(m-4) \otimes \mathrm{MU}(2) \otimes \mathrm{MU}(2) \\ &= \mathrm{MTSU}(m-4) \otimes \Sigma^6 B(U(2) \times U(2))_+^{\overline{V}_{U(2)} \oplus \overline{V}_{U(2)}} \\ &\xrightarrow[(8.3)]{\widetilde{\mathbb{CP}}_{2,2}\cdot} \mathrm{MTSU}(m-4) \otimes \mathrm{MTSU}(4) \\ &\xrightarrow{\mathrm{mult}} \mathrm{MTSU}(m)\end{aligned}$$

For $m < 4$, we define $\mathrm{flop}_m \colon * \to \mathrm{MTSU}(m)$ to be the trivial map. We call the cofiber

$$\mathrm{MTSU}(m)/\mathrm{flop}_m$$

the *tangential* $\mathrm{SU}(m)$*-bordism spectrum modulo flops*. The map $\mathrm{flop}_m$ is compatible with the stabilization $m \mapsto m+1$, so we define

$$\mathrm{MSU}/\mathrm{flop} = \operatorname*{hocolim}_{m \to \infty} \mathrm{MTSU}(m)/\mathrm{flop}_m$$

and call it *the stable* SU*-bordism spectrum modulo flops*.

*Remark* 8.5*.* Definition 8.4 agrees with the notion of SU-flops in [Tot00]. Indeed, $\mathrm{flop}_m$ is geometrically understood as follows. An element in $\pi_{n-6} \, \mathrm{MTSU}(m-4) \otimes \mathrm{MU}(2) \otimes \mathrm{MU}(2)$ is represented by a smooth $(n-6)$-dimensional manifold $Z$ equipped with two rank-two complex vector bundles $V_A$ and $V_B$ together with a stable $SU(m-4)$-structure on $TZ \oplus V_A \oplus V_B$. The map $\mathrm{flop}_m$ sends it to the manifold $\widetilde{\mathbb{CP}}(V_A \oplus V_B)$, the $n$-dimensional manifold with stable tangential $SU(m)$-structure, obtained by performing the bundle version of the construction in Definition 8.1. This is exactly the procedure described in [Tot00, pp. 773, 777]. Totaro defines the ideal $I \subset \pi_* \, \mathrm{MSU}$ generated by manifolds of the form $\widetilde{\mathbb{CP}}(V_A \oplus V_B)$, and call $\pi_* \, \mathrm{MSU}/I$ the SU-bordism group modulo flops. ⌟

8.2. **Invariance under** SU**-flops.** We now show that the $U(1)$-equivariant topological elliptic genera with values in (periodic) $\mathrm{TJF}_m$ is SU-flop invariant:

**Theorem 8.6.** *The topological elliptic genera factor, naturally in $m$, as maps*

$$\mathrm{Jac}_m \colon \, \mathrm{MTSU}(m)/\mathrm{flop}_m \to \mathrm{TJF}_m,$$

*In particular, we get a stable factorization*

$$\mathrm{Jac} \colon \, \mathrm{MSU}/\mathrm{flop} \to \mathrm{TJF}.$$

The key observation for constructing such a factorization is the following. Let us use the notation

$$\mathrm{TJF}_m[V]^G := \mathrm{TJF}_m \otimes_{\mathrm{TMF}} \mathrm{TMF}[V]^G$$

for $V \in \mathrm{RO}(G)$.

**Lemma 8.7.** *We have*

$$\mathrm{TMF}[-V_{U(2)}]^{U(2)} = 0$$

*as a* TMF*-module spectrum. In particular, we have*

$$\mathrm{TJF}_4 \left[- \left(\overline{V}_{U(2)} \oplus \overline{V}_{U(2)}\right)\right]^{U(2) \times U(2)} = 0,$$

*so that we have*

$$\mathrm{Jac}_4[\widetilde{\mathbb{CP}}_{2,2}] = 0 \in \pi_6 \mathrm{TJF}_4 \left[-\left(\overline{V}_{U(2)} \oplus \overline{V}_{U(2)}\right)\right]^{U(2)\times U(2)}.$$

*Proof of Lemma 8.7.* We use the stabilization-restriction fiber sequence of [LY24, Section 4],

$$\mathrm{TMF}[-V_{U(2)}]^{U(2)} \xrightarrow{\chi(V_{U(2)})\cdot} \mathrm{TMF}^{U(2)} \xrightarrow{\mathrm{res}^{U(1)}_{U(2)}} \mathrm{TMF}^{U(1)}.$$

for the restriction along the inclusion $\iota\colon U(1) \hookrightarrow U(2)$ to one of the diagonal entry. It is enough to show that the second arrow is an equivalence. Since $\iota$ has a retraction by the determinant, the following composition is the identity on $\mathrm{TMF}^{U(1)}$:

$$\mathrm{id}\colon \mathrm{TMF}^{U(1)} \xrightarrow{\mathrm{res}_{\det}} \mathrm{TMF}^{U(2)} \xrightarrow{\mathrm{res}^{U(1)}_{U(2)}} \mathrm{TMF}^{U(1)}. \tag{8.8}$$

By [GM] (see also [LY24, Fact 6.5 (2)]), the first arrow in (8.8) is an isomorphism, thus so is the second arrow. □

*Proof of Theorem 8.6.* We construct the factorization as follows. First we construct it in the case $m = 4$. For that, we use the following commutative diagram of fiber sequences

$$\begin{array}{ccccc}
\Sigma^6 \mathrm{MU}(2) \otimes \mathrm{MU}(2) & \xrightarrow{\mathrm{flop}_4} & \mathrm{MTSU}(4) & \longrightarrow & \mathrm{MTSU}(4)/\mathrm{flop} \\
\downarrow{\scriptstyle \mathrm{Jac}_0} & & \downarrow{\scriptstyle \mathrm{Jac}_4} & & \downarrow{\scriptstyle \mathrm{Jac}_4} \\
\Sigma^6 \mathrm{TJF}_0[\overline{V}_{U(2)} \boxplus \overline{V}_{U(2)} - \mathrm{ad}]^{U(2)\times U(2)} & \xrightarrow{\mathrm{Jac}(\widetilde{\mathbb{CP}}_{2,2})\cdot} & \mathrm{TJF}_4 & \longrightarrow & \mathrm{TJF}_4/\mathrm{Jac}(\widetilde{\mathbb{CP}}_{2,2})
\end{array}$$

The left vertical arrow is given by tensoring $\mathrm{Jac}_0\colon S = \mathrm{MTSU}(0) \to \mathrm{TJF}_0$ and the norm map $BG^V \to \mathrm{TMF}[V - \mathrm{ad}(G)]^G$. However, we have

$$\begin{aligned}
&\mathrm{Map}_{\mathrm{Mod}_{\mathrm{TJF}_0}}\left(\Sigma^6 \mathrm{TJF}_0[\overline{V}_{U(2)} \boxplus \overline{V}_{U(2)} - \mathrm{ad}]^{U(2)\times U(2)}, \mathrm{TJF}_4\right) \\
&\simeq \Sigma^{-6}\mathrm{TJF}_4\left[-\left(\overline{V}_{U(2)} \oplus \overline{V}_{U(2)}\right)\right]^{U(2)\times U(2)} = 0.
\end{aligned}$$

where the last equality is Lemma 8.7. This gives the splitting $\mathrm{TJF}_4/\mathrm{Jac}(\widetilde{\mathbb{CP}}_{2,2}) \to \mathrm{TJF}_4$ which gives the desired factorization for $m = 4$.

For the general case $m \geq 4$, we use the following diagram,

$$\begin{array}{ccccc}
\mathrm{MTSU}(m-4) \otimes \mathrm{MTSU}(4) & \to & \mathrm{MTSU}(m-4) \otimes \mathrm{MTSU}(4)/\mathrm{flop}_4 & \xrightarrow{\mathrm{Jac}_{m-4}\otimes \mathrm{Jac}_4} & \mathrm{TJF}_{m-4} \otimes \mathrm{TJF}_4 \\
\downarrow & & \downarrow & & \\
\mathrm{MTSU}(m) & \longrightarrow & \mathrm{MTSU}(m)/\mathrm{flop}_m & & \\
\downarrow{\scriptstyle \mathrm{Jac}_m} & & \downarrow{\scriptstyle \mathrm{Jac}_m}\ (\text{dotted}) & \nwarrow{\scriptstyle \mathrm{multi}} & \\
\mathrm{TJF}_m & = & \mathrm{TJF}_m & &
\end{array}$$

where the top horizontal arrows use the case $m = 4$ above, and solid parts of the diagram commute. The left upper square is pushout, so the dotted arrow is induced, which gives the desired factorization $\mathrm{Jac}_m$. The naturality in $m$ is obvious by construction.

□

KTH Royal Institute of Technology, Stockholm, Sweden
*Email address*: `tilmanb@kth.se`

Perimeter Institute for Theoretical Physics, 31 Caroline Street North, Waterloo, Ontario, Canada, N2L 2Y5
*Email address*: `myamashita@perimeterinstitute.ca`